\documentclass[12pt]{amsart}
\usepackage{amsmath}\usepackage{amsrefs}
\usepackage{amsfonts}\usepackage{dictsym}
\usepackage{amssymb}\usepackage{marvosym}
\usepackage{verbatim}\usepackage{MnSymbol}
\usepackage{skak}%\usepackage{hyperref}
\usepackage{dsfont}
\usepackage{amstext}
\usepackage{amsbsy}
\usepackage{amsopn}
\usepackage{amsthm}\usepackage{color}

\theoremstyle{definition}

\theoremstyle{remark}

\def\proclaim#1{\vskip0.5em\noindent{\bf #1}\it}
\def\endproclaim{\vskip0.5em\par\noindent\rm}
\def\stackunder#1#2{\mathrel{\mathop{#2}\limits_{#1}}}%
\def\proclaim#1{\vskip0.5em\noindent{\bf #1}\it}
\def\endproclaim{\vskip0.5em\par\noindent\rm}
\def\stackunder#1#2{\mathrel{\mathop{#2}\limits_{#1}}}%
\def\demo#1{\vskip0.5em\noindent{\bf #1\ }}

\def\text#1{\mbox{#1}}
\def\flushpar{\par\noindent}

\def\tag#1{\eqno{(#1)}}
\def\mod{\mbox{ mod }}
\newcommand{\mapright}[1]{%
    \smash{\mathop{%
        \hbox to 1cm{\rightarrowfill}
        }
    \limits^{#1}
    }
}
\newcommand{\mapleft}[1]{%
    \smash{\mathop{%
        \hbox to 1cm{\rightarrowfill}
        }
    \limits_{#1}
    }
}

\usepackage{wasysym}\usepackage{eiad}%\usepackage{sauter}
\def\e{\epsilon}
\def\a{\alpha}
\def\b{\beta}
\def\g{\gamma}

\def\d{\delta}
\def\D{\Delta}
\def\th{\theta}
\def\s{\sigma}

\def\l{\lambda}

\def\x{\times}

\def \C{C \! \! \! \! I \ }
\def\o{\overline}
\def\f{\flushpar}
\def\u{\underline}
\def\v{\varphi}

\def\bdy{\partial}

\def\om{\omega}
\def\Om{\Omega}
\def\B{\Cal B}
\def\C{\Cal C}

\def\({\biggl(}
\def\){\biggr)}
\def\[{\bigl[}
\def\]{\bigr]}\def\Cal{\mathcal}

\def\<{\langle}
\def\>{\rangle}

\def\wrt{\text{with respect to}}

\def\bul{\smallskip\f$\bullet\ \ \ $}\def\Par{\smallskip\f\P}
\def\pf{\f{\it Proof}\ \ \ \ }\def\sms{\smallskip\f}\def\sbul{\f$\bullet\
\ \ $}\def\sms{\smallskip\f}\def\wrt{\text{with respect to}}
\def\st{\text{such that}}\def\lra{\longrightarrow}

\def\Lra{\Longrightarrow}

\def\lfl{\lfloor}\def\rfl{\rfloor}

\begin{document}
 \title{ Rational weak mixing in infinite measure spaces}
\author{ Jon. Aaronson}
 \address{\ \ School of Math. Sciences, Tel Aviv University,
69978 Tel Aviv, Israel.\f\ \ \  {\it Webpage }: {\tt http://www.math.tau.ac.il/$\sim$aaro}}\email{aaro@post.tau.ac.il}
\thanks{Research
supported by Israel Science Foundation grant No. 1114/08. }
\thanks {\copyright Jon Aaronson 2010-11}

\markboth{\copyright Jon Aaronson 2010-11}{Rational weak mixing}
\begin{abstract}
 Rational weak mixing is a  measure theoretic version of Krickeberg's   strong ratio mixing property for infinite measure preserving transformations. It requires ``{\tt density}" ratio convergence for every pair of  measurable sets in a dense hereditary
ring. Rational weak mixing implies weak rational ergodicity and (spectral) weak mixing. It is enjoyed for example  by
Markov shifts with  Orey's    strong ratio limit property. The power, subsequence version of the property is generic.

\end{abstract}\maketitle%\end{document}

\section*{\S0  Introduction:\ \ Hopf's example}

E. Hopf gave an example in [H] of a transformation of the infinite strip $\mathbb R_+\x [0,1]$, preserving Lebesgue measure $m$ which satisfies the {\tt ratio mixing} property:
\begin{align*}\tag{$\maltese$}\label{maltese}&\frac{m(A\cap T^{-n}B)}{u_n}\underset{n\to\infty}\lra m(A)m(B)\\ & \ \forall\ \ A,\ B\ \text{\tt bounded with}\ \ m(\bdy A)=m(\bdy B)=0\end{align*}
where $u_n=\sqrt{\frac2{\pi n}}$.

Hopf mentioned that if ($\maltese$)   could be established for every bounded measurable set, this would imply ergodicity of $T$. This latter property is also invariant under isomorphism.
 \

 The theory of {\tt weakly wandering sets} as in [HK] shows that ($\maltese$)  cannot hold for every pair of sets in any dense,
{\tt hereditary collection} (see below) in the absence of absolutely continuous, invariant probabilities
(and  this cannot be used to establish ergodicity of Hopf's example).

 \

 We show here that Hopf's example is {\tt rationally weakly mixing} in a sense which implies that for every pair of bounded measurable sets, ($\maltese$) (as on page  \pageref{maltese}) is satisfied on a subsequence of full density.
 \

\subsection*{Organization of the paper}
\

In \S1, we give  definitions and  preliminary discussions. The main results are stated in \S2.
In \S3, we study the modes of convergence involved in the rational weak mixing properties. \S4 contains the proof of the ``basic" proposition 0 and the
``density convergence" theorem A.
 In \S5, we establish sufficient conditions for rational weak mixing (Lemmas B and C). We  collect in \S6 some facts on the mean ergodic theorem with
 weighted averages for use in \S7 to prove theorem D which connects subsequence rational weak mixing with other
mixing properties. Markov shifts are
treated in \S8 (where there is some discussion of smoothness of renewal sequences) and Gibbs-Markov towers are studied in \S9 via their local limit properties.
We ``make categorical statements'' (theorem F)  in \S10 and some closing remarks in \S11.
\

 \section*{\S1 Definitions and preliminaries}
 \subsection*{Notation and basics}
 \

 In this paper $(X,\B,m,T)$ denotes a measure preserving transformation $T$ of a non-atomic, $\s$-finite, standard measure space $(X,\B,m)$.

  Unless otherwise stated, the measure will be infinite ($m(X)=\infty$). Measure preserving transformations of finite measure spaces are referred to as probability preserving transformations.

  A standing assumption on $(X,\B,m,T)$ is {\it conservativity}: $$m(A\setminus\bigcup_{n=1}^\infty T^{-n}A)=0\ \forall\ A\in\B.$$
  The collection of measurable subsets of $X$ with finite measure is denoted  $\mathcal F:=\{A\in\B:\ m(A)<\infty\}$ and for any $\mathcal C\subset\B$ the ``positive elements" of $\mathcal C$ are denoted
 $$\mathcal C_+:=\{A\in\mathcal C:\ m(A)>0\}.$$

\subsection*{Krickeberg mixing}\ \

Krickeberg ([Kri1]) noted that Hopf's example is isomorphic to the Markov shift of the simple symmetric random walk on $\mathbb N$ with reflecting barrier at $1$ which has irreducible, recurrent transition matrix ([KM]) and is therefore conservative, ergodic ([HR]). He also formulated a concept of topological ratio mixing for transformations preserving infinite measures:
\

Let $(X,\B,m,T)$ be a measure preserving transformation and let $\a\subset\B$ be a countable partition, generating $\B$ under $T$ in the sense that $\s(\bigcup_{n\in\mathbb Z}T^n\a)=\B$. The measure preserving transformation $(X,\B,m,T)$ is called {\it Krickeberg $\a$-mixing}  if \sms $\exists\ u_n>0\ \ (n\ge 1)\ \st$ ($\maltese$) (as on page  \pageref{maltese})  is satisfied $\forall\ A,\ B\in\mathcal C_\a$, the collection of $(\a,T)$-{\it cylinder sets} defined by
$$\mathcal C_\a=\mathcal C_\a(T):=\{[a_1,\dots,a_N]_k:=\bigcap_{j=1}^NT^{j_k}a_j:\ N\in\mathbb N,\ k\in\mathbb Z,\ a_1,\dots,a_N\in\a\};$$
and hence ([Kri1])  $\forall\ \ A,\ B$  with $m(\bdy A)=m(\bdy B)=0$ when $X$ is considered equipped with the product topology from $\a^\mathbb Z$.
\

Markov shifts with the {\it strong ratio limit property} (SRLP) as in [O] (e.g. Hopf's example) are Krickeberg $\a$-mixing with
$\a$ the natural partition according to the state occupied at time $0$ ([Kri1]). Examples of
Krickeberg $\a$-mixing measure preserving transformations are also given in [Fri], [Pap], [T] and [MT]. Other definitions of mixing are discussed in [L].
\

It follows from theorem 8.1  that Markov shifts whose associated renewal sequences have the   strong ratio limit property   are {\tt rationally weakly mixing}.

\subsection*{Hereditary rings}\ \
\

Let $(X,\B,m)$ be a $\s$-finite measure space. A collection $\mathcal C\subset\B$ is called    {\it hereditary} if
 $$\mathcal H(\C):=\{A\in\B,\ A\subseteq B\in\mathcal C\}=\C.$$ A {\it hereditary ring} $\mathcal H\subset\B$ is a hereditary collection which is closed under finite union. It is {\it dense} if
$$\forall\ A\in\mathcal F :=\{F\in\B:\  m(F)<\infty\},\ \e>0\ \exists\ H\in\mathcal H,\ m(A\setminus H)<\e.$$
 \

For example, both
 $\mathcal F $ and the collection $\mathcal R_b$ of bounded measurable subsets of the infinite strip $\mathbb R_+\x [0,1]$ are dense hereditary rings. The collection of null sets is a hereditary ring which is not dense. We'll denote the minimal hereditary ring containing the collection $\mathcal C\subset\B$ by
 $\mathcal{HR(C)}$.
 \

Any two dense, hereditary rings in the same measure space intersect and thus many ergodic properties
involving such are isomorphism invariant (e.g. rational weak mixing).

\subsection*{Standard measure spaces}
\

We assume that all measure preserving transformations are defined on  standard measure spaces. The $\s$-finite  measure space $(X,\B,m)$ is {\it standard} if $X$ is a Polish space, $\B$ is the collection of Borel sets and $m$ is non-atomic. The standardness assumption is used as follows:
\bul If $(X,\B,m,T)$ is a conservative, ergodic, measure preserving transformations of a standard measure space then $\exists$ a countable partition $\a\subset\mathcal F:=\{A\in\B:\ m(A)<\infty\}$ which generates $\B$ under $T$ and, up to isomorphism, $T$ is the shift on  $X=\a^\mathbb Z$ equipped with the product topology (a homeomorphism).
\subsection*{Weights}
\

Our results involve averaging techniques using certain non-negative, bounded  {\tt weight} sequences. We call a bounded sequence
$u=(u_0,u_1,\dots)$ an {\it admissible weight sequence} (abbr. to {\it weight}) if
$$u_n\ge 0\ \forall\ n\ge 1\ \&\ \  a_u(n):=\sum_{k=1}^nu_k\lra \infty$$
and denote the collection of weights by $\mathfrak W.$
\

We'll denote, for (eventually) positive sequences $u=(u_0,u_1,\dots)\ \&\ w=(w_0,w_1,\dots)$:
\bul $u_n\sim w_n$ if $\tfrac{u_n}{w_n}\underset{n\to\infty}\lra\ 1$;

and for non-negative sequences $u=(u_0,u_1,\dots)\ \&\ w=(w_0,w_1,\dots)$:
\bul $u_n\ll w_n$ if $\exists\ M>0\ \st\ u_n\le Mw_n\ \forall\ n\ge 0$;
\bul $u_n\asymp w_n$ if  $u_n\ll w_n$ and $u_n\gg w_n$.

\

Given a subsequence $\mathfrak K\subset\mathbb N$, we call weights $u,\ w\in\mathfrak W$ \sbul  {\it $\mathfrak K$-asymptotic  }  ($u\overset{\mathfrak K}\approx w$) if
$\frac1{a_u(n)}\sum_{k=1}^n|u_k-w_k|\underset{n\to\infty,\ n\in\mathfrak K}\lra 0.$
\

Evidently, $$u\overset{\mathfrak K}\approx w\ \implies\ \frac{a_w(n)}{ a_u(n)}\underset{n\to\infty,\ n\in\mathfrak K}\lra 1.$$
 The converse implication sometimes holds and will be discussed in the sequel.
\

We call $u\in\mathfrak W$  \bul {\it $\mathfrak K$-smooth  }   if  $\frac1{a_u(n)}\sum_{k=1}^n|u_k-u_{k+1}|\underset{n\to\infty,\ n\in\mathfrak K}\lra 0$; equivalently
$(u_0,u_1,\dots)\overset{\mathfrak K}\approx (u_1,u_2,\dots)$.

\

We'll say that weights $u,\ w\in\mathfrak W$  are {\it asymptotic} ($u\approx w$) if $u\overset{\mathbb N}\approx w$, that   $u\in\mathfrak W$ is {\it smooth}
if it is   $\mathbb N$-smooth and {\it subsequence smooth} if it is $\mathfrak K$-smooth for some subsequence $\mathfrak K\subset\mathbb N$.
\

\subsection*{Intrinsic weights}\ \

For $(X,\B,m,T)$   a  conservative, ergodic, measure preserving transformation and  $E,\ F\in\mathcal F_+$ the {\it intrinsic weight}  $u(E,F)\in\mathfrak W$ is defined by
$$u_n(E,F)\ :=\ \frac{m(F\cap T^{-n}F)}{m(E)m(F)}.$$We denote $a_n(E,F)=a_{u(E,F)}(n):=\tsum_{k=0}^{n-1}u_k(E,F)$ and write $u(F):=u(F,F)$.

\subsection*{Rational weak mixing}\ \

Let $\mathfrak K\subset\mathbb N$ be a subsequence. We'll call the conservative, ergodic, measure preserving transformation $(X,\B,m,T)$   {\it rationally weakly mixing along $\mathfrak K$}  if
$\exists\ F\in\mathcal F_+$ so that
\begin{align*}\tag{$\bigstar_{\mathfrak K}$}\label{bigstarK}m(A\cap T^{-n}B)\ \ \overset{\mathfrak K}\approx\ \ \ m(A)m(B)u_n(F)\ \ \forall\ A,\ B\in \B\cap F.\end{align*}
We call the measure preserving transformation $(X,\B,m,T)$ \bul {\it rationally weakly mixing}  if it is    rationally weakly mixing along $\mathbb N$; and
\bul   {\it subsequence rationally weakly mixing}   if it is    rationally weakly mixing along some $\mathfrak K\subset\mathbb N$.

\subsection*{Weak rational ergodicity}\ \
\

Again for $\mathfrak K\subset\mathbb N$ a subsequence,
the conservative, ergodic, measure preserving transformation $(X,\B,m,T)$ is called {\it weakly rationally ergodic along $\mathfrak K$}  if $\exists$   $F\in\mathcal F_+$ so  that
\begin{align*}\tag{$\largestar_\mathfrak K$}\label{largestarK}\frac1{a_n(F)}\sum_{k=0}^{n-1}m(B\cap T^{-k}C)\underset{n\to\infty,\ n\in\mathfrak K}\lra\ m(B)m(C)\ \forall\ B,C\in\B\cap F\end{align*}
where $a_n(F):=\frac1{m(F)^2}\sum_{k=0}^{n-1}m(F\cap T^{-k}F)$. Weak rational ergodicity entails conservativity and ergodicity.
\

The proof of theorem 3.3 in [FL] easily adapts to show that $F\in\mathcal F_+$ satisfies ($\largestar_\mathfrak K$) if and only if
\begin{align*}\{\tfrac{S_n^{(T)}(1_F)}{a_n(F)}:\ n\in\frak K\}\ \ \text{\tt is uniformly integrable on}\ F.\end{align*}
A useful sufficient condition for this is
$$\sup_{n\in\frak K}\tfrac1{a_n(F)^2}\int_FS_n(1_F)^2dm<\infty$$ and  $(X,\B,m,T)$ is called {\it rationally ergodic along $\mathfrak K$}  if $\exists$   $F\in\mathcal F_+$ with this property. See
 [A], [A1] (for the special case $\frak K=\mathbb N$).
\

In case $T$ is weakly rationally ergodic along $\mathfrak K$:
\sbul the collection of sets $R_{\mathfrak K}(T)$ satisfying ($\largestar_\mathfrak K$)  is  a hereditary ring;
\sbul  $\exists\ a_n(T)$ (the {\it return sequence along $\mathfrak K$}) $\st$
$$\frac{a_n(A)}{a_n(T)}\underset{n\to\infty,\ n\in\mathfrak K}\lra 1\ \forall\ A\in R_{\mathfrak K}(T);$$
\sbul for conservative, ergodic $T$,  $R_{\mathfrak K}(T)=\mathcal F$ only when $m(X)<\infty$.
\
The proofs of these statements  are analogous to those
 in [A], [A1] (for the special case $\frak K=\mathbb N$).
 \

\

We'll call a measure preserving transformation $(X,\B,m,T)$ :
\bul {\it [weakly] rationally ergodic} if it is [weakly] rationally ergodic along $\mathbb N$ and set $R(T):=R_{\mathbb N}(T)$ (as in [A], [A1]); and
\bul {\it subsequence [weakly] rationally ergodic} if it is [weakly] rationally ergodic along some $\mathfrak K\subset\mathbb N$.

For example,  conservative, ergodic Markov shifts are  rationally ergodic. For further examples,  see [A].

\

We'll see that   rational weak mixing along $\mathfrak K$ implies weak rational ergodicity  along $\mathfrak K$  and   that for $T$ rationally weakly mixing along $\mathfrak K$  ,
$$\{F\in\mathcal F_+:\ \ \text{\tt($\bigstar_{\mathfrak K}$) holds}\}=R_{\mathfrak K}(T)$$ where ($\bigstar_{\mathfrak K}$) is as on page  \pageref{bigstarK}.

\

\subsection*{Weak mixing}
\

For  a measure preserving transformation $(X,\B,m,T)$ of a $\s$-finite measure space (as shown in [ALW])
the following conditions are equivalent:
\begin{align*}\tag{i}f\in L^\infty,\ \l\in\mathbb S^1,\ f\circ T=\l f\ \text{\tt a.e.}\ \Rightarrow\ f\ \text{\tt is constant a.e.}\end{align*}
\begin{align*}\tag{ii} T\x S\ \text{\tt\small is ergodic}\ \forall\ \text{\tt\small  ergodic, probability preserving}\ S;\end{align*}
\begin{align*}\tag{iii}\tfrac1n\sum_{k=0}^{n-1}|\int_Xuf\circ T^kdm|\underset{n\to\infty}\lra 0\ \forall\ u\in L^1_0,\ f\in L^\infty.\end{align*}

We'll call a  measure preserving transformation  satisfying (any one of) them {\it spectrally weakly mixing}. This in the interest of disambiguation.
In [ALW] and elsewhere
``spectral weak mixing" is called  "weak mixing".
\

We'll see that subsequence rational weak mixing $\implies$ spectral weak mixing.
\subsection*{Categorical statements}
\

Let $(X,\B,m)$ be a standard $\s$-finite, non-atomic, infinite measure space and consider $\text{\tt MPT}(X,\B,m)$,  the collection of  invertible measure preserving
transformations of $(X,\B,m)$
equipped with the {\it weak operator topology } defined by $T_n\to T$ if
$$m(TA\D T_nA)+m(T^{-1}A\D T_n^{-1}A)\underset{n\to\infty}\lra\ 0\ \ \forall\ A\in\mathcal F.$$  It follows that  $\text{\tt MPT}(X,\B,m)$ is a Polish group.
A {\it categorical statement} is a statement concerning the Baire category of a subset of {\tt MPT}. For a review of this subject, see [CP].
\

 We'll see that (the {\tt power} version of)
the subsequence rational weak mixing elements of {\tt MPT} form a  residual set in {\tt MPT}.

\section*{\S2 Results}

\proclaim{Proposition 0\ \ ({\tt basics})}
\

Let $\mathfrak K\subset\mathbb N$ be a subsequence and suppose that $(X,\B,m,T)$ is rationally weakly mixing along $\mathfrak K$, then

\sms {\rm (i)}\ \ $T$ is weakly rationally ergodic along $\mathfrak K$;
\begin{align*}\tag{ii} \{F\in\mathcal F_+:\ \ \text{\tt($\bigstar_{\mathfrak K}$) holds}\}=R_{\mathfrak K}(T)\end{align*}
 where ($\bigstar_{\mathfrak K}$) is as on page  \pageref{bigstarK};
\sms {\rm (iii)}\ \ \   $u(F,G)$ is $\mathfrak K$-smooth   $\forall\ F,\ G\in R_{\mathfrak K}(T),\ m(F),\ m(G)>0$.
\sms {\rm (iv)}\ \ \  for each $p\in\Bbb Z,\ p\ne 0,\ T^p$ is rationally weakly mixing along  $\tfrac1{|p|}\mathfrak K:=\{\lfl\tfrac{k}{|p|}\rfl:\ k\in\frak K\}$ and $R_{\tfrac1{|p|}\mathfrak K}(T^p)=R_{\mathfrak K}(T)$.
\endproclaim

.
\subsection*{Corollary: Isomorphism invariance}
\

It follows from Proposition 0 that if $(X,\B,m,T)$ is rationally weakly mixing along $\mathfrak K$,  and is isomorphic by measure preserving transformation to $(X',\B',m',T')$,
then $(X',\B',m',T')$ is also  rationally weakly mixing along $\mathfrak K$,
and the intrinsic weights of $T'$ are $\mathfrak K$-asymptotic to those of $T$.
\proclaim{Theorem A\ \ ({\tt density convergence})}
\

Suppose that $(X,\B,m,T)$ is rationally weakly mixing and that
\sbul $\exists\ E\in R_{\mathfrak K}(T)$ st  $u(E)\approx v$ where $v\in\mathfrak W$ is {\it regularly varying} with {\it index} $s\in (-1,0]$
($\tfrac{v_{\lfl\l n\rfl}}{v_n}\underset{n\to\infty}\lra\l^s\ \ \forall\ \l>0$),  then for $F\in R_{\mathfrak K}(T)_+$,
\begin{align*} \tfrac{m(A\cap T^{-n}B)}{u_n(F)}\overset{\text{\tt\tiny density}}{\underset{n\to\infty}\lra}\ m(A)m(B)\ \ \forall\ A,\ B\in R_{\mathfrak K}(T).\end{align*}
\endproclaim
Here  $s_n\overset{\text{\tt\tiny density}}{\underset{n\to\infty}\lra}L$ means $s_n{\underset{n\to\infty,\ n\notin K}\lra}L$ where $K\subset\mathbb N$ has {\it zero density} in  the sense that $\#(K\cap [1,n])=o(n)$ as $n\to\infty$.

\

\proclaim{Lemma B\ \ ({\tt sufficient conditions})}\ \ Let $\mathfrak K\subset\mathbb N$ be a subsequence. Suppose that $X$ is a Polish space and
\sms {\rm (i)} $(X,\B,m,T)$ is an invertible,  measure preserving transformation, weakly rationally ergodic  along $\mathfrak K$;
\sms {\rm (ii)} $\exists\ \Om\in R_\frak K(T)$ open in $X$ and a countable base $\mathcal C$ for the topology of $\Om\ \st$

{\small\begin{align*}\tag{a}\forall\ \{C_i\}_{i=1}^k\subset\mathcal C\ \exists\ \{D_j\}_{j=1}^\ell\subset\mathcal C,\ m(D_i\cap D_j)=0\ \forall\ i\ne j;
\bigcup_{i=1}^kC_i\overset{m}=\bigcup_{j=1}^\ell D_j;\end{align*} }
\begin{align*}\tag{b}\ m(A\cap T^{-n}B)\ \ \overset{\mathfrak K}\approx\ \ \ m(A)m(B)u_n(\Om)\ \ \forall\ A,\ B\in\mathcal C.\end{align*}
 then $(X,\B,m,T)$ is rationally weakly mixing along $\mathfrak K$.
\endproclaim
\

\proclaim{Lemma C\ \ ({\tt sufficient conditions})}\ \ Let $\mathfrak K\subset\mathbb N$ be a subsequence. Suppose that
\sms {\rm (i)} $(X,\B,m,T)$ is an invertible,  measure preserving transformation, weakly rationally ergodic  along $\mathfrak K$;
\sms {\rm (ii)} $\exists$ a countable generating partition $\a\subset R_{\mathfrak K}(T)$ and $\Om\in\mathcal C_\a\ \st$
\begin{align*}m(A\cap T^{-n}B)\ \ \overset{\mathfrak K}\approx\ \ \ m(A)m(B)u_n\ \ \forall\ A,\ B\in\mathcal C_\a\end{align*} where
$u=u(\Om)$,
 then $(X,\B,m,T)$ is rationally weakly mixing along $\mathfrak K$.
\endproclaim
\

\

\proclaim{Theorem D\ \ ({\tt mixing properties})}
\

Let $\mathfrak K\subset\mathbb N$ be a subsequence and suppose that $(X,\B,m,T)$ is rationally weakly mixing along $\mathfrak K$, then
\sms {\rm (i)} $(X,\B,m,T)$ is spectrally weakly mixing;
\sms {\rm (ii)} $T\x S$ is rationally weakly mixing  along $\mathfrak K$ $\forall$ weakly mixing, probability preserving transformation $(\Om,\mathcal F,P,S)$.\endproclaim

Invertible rationally weakly mixing measure preserving transformations of infinite measure spaces are obtained via

\proclaim{Corollary E}
\

 The natural extension of a measure preserving transformation, rationally weakly mixing along $\mathfrak K$,  is also rationally weakly mixing along $\mathfrak K$ with  $\mathfrak K$-asymptotic
 intrinsic weights.
\endproclaim

\

Let
\begin{align*}&\text{\tt RWM}(X):=\{T\in\text{\tt MPT}(X):\ T\ \text{\tt\small is rationally weakly mixing}\}\\ &\,\\ & \text{\tt SRWM}(X):=\\ &\{T\in\text{\tt MPT}(X):\ T\ \text{\tt\small is subsequence rationally weakly mixing}\}.\end{align*}

\

For $T\in\text{\tt MPT}(X),\ \D\ge 1$ and $(\kappa_1,\dots,\kappa_\D)\in\Bbb Z^\D,$ let
$$T^{(\kappa_1,\dots,\kappa_\D)}:=T^{\kappa_1}\x T^{\kappa_2}\x\dots\x T^{\kappa_\D} \in\text{\tt MPT}(X^\D).$$
Call $T\in\text{\tt MPT}$ {\it power, subsequence rationally weakly mixing} if $T^{(\kappa_1,\dots,\kappa_\D)}$ is subsequence rationally weakly mixing
\

$\forall\ \D\ge 1\ \&\ (\kappa_1,\dots,\kappa_\D)\in(\Bbb Z\setminus\{0\})^\D.$ Let
\begin{align*}&\text{\tt PSRWM}(X):=\\ &\{T\in\text{\tt MPT}(X):\ T^{(\kappa_1,\dots,\kappa_\D)}\in\text{\tt SRWM}(X^\D)\ \forall\ \D\ge 1,\ (\kappa_1,\dots,\kappa_\D)\in(\Bbb Z\setminus\{0\})^\D\}.\end{align*}
\proclaim{Theorem F\ \ ({\tt Baire category})}
\

\sms {\rm (i)}\ \ The collection {\tt RWM} is meagre in $\text{\tt MPT}$.

\sms {\rm (ii)}\ \ The collection
$\text{\tt PSRWM}(X)$ is residual in $\text{\tt MPT}(X)$.
\endproclaim
\section*{\S3 Convergence}
\

In this section we study the modes of convergence involved in the rational weak mixing properties. Let $\mathfrak K\subset\mathbb N$ be a subsequence.
\

\subsection*{$(u,\mathfrak K)$-small sets}
\

Let $u\in\mathfrak W$.
We'll say that the set $K\subset \mathbb N$
\bul is {\it $(u,\mathfrak K)$-small} if $\frac{a_u(K,n)}{ a_u(n)}\underset{n\to\infty,\ n\in\mathfrak K}\lra 0$ where
$ a_u(K,n):=\sum_{k\in K\cap [1,n]}u_k$; and  {\it $(u,\mathfrak K)$-large } if $K^c$ is {\it $(u,\mathfrak K)$-small}.

\

We'll call a set {\it $u$-small} if it is $(u,\mathbb N)$-small.
\

Recall that the set $K\subset \mathbb N$\bul  has {\it density} $d(K)$ if $\tfrac1n\#(K\cap [1,n])\lra d(K)$ as $n\to\infty$; \bul and has {\it zero density} if $d(K)=0$ (equivalently: $K$ is $\mathds{1}$-small where $\mathds{1}\in \mathfrak W,\ \mathds{1}_n= 1\ \forall\ n$).

\

The following remark collects some elementary facts about $(u,\mathfrak K)$-smallness:

\subsection*{Remark 3.1}
\

Suppose that $u\in \mathfrak W$, then
\sms {\rm(o)}  if $K\subset \mathbb N$
 is  $(u,\mathfrak K)$-small, it is also $(u,\mathfrak K')$-small whenever $\mathfrak K'\subset\mathfrak K+F$ for some finite set $F\subset\Bbb Z$.
This is because if $n\in\mathfrak K'\ \&\ n=m+f$ where $m\in\frak K\ \&\ f\in f$, then
 $$|a_u(K,n)-a_u(K,m)|\le \sup_k u_k\cdot\max_{j\in F}|j|.$$
\sms {\rm (i)} a finite union of $(u,\mathfrak K)$-small sets is itself $(u,\mathfrak K)$-small;

\sms {\rm (ii)}  if $v\in\mathfrak W$ satisfies $u_n\asymp v_n$ then $K\subset \mathbb N$ is $(u,\mathfrak K)$-small iff it is $(v,\mathfrak K)$-small;
\sms {\rm (iii)}\label{rem3.1.3}  if $K_1\subset K_2\subset\dots$ is an increasing sequence of $(u,\mathfrak K)$-small sets, then $\exists\ N_1<N_2<\dots$ so that $K_\infty:=\bigcup_{j=1}^\infty K_j\cap [N_j+1,N_{j+1}]$ is a $(u,\mathfrak K)$-small set.
\pf A suitable sequence is obtained by  choosing $N_j\uparrow$ $\st$
$\tfrac{ a_u(K_j,n)}{ a_u(n)}<\tfrac1j\ \forall\ n\ge N_j,\ n\in\mathfrak K$.
\sms {\rm (iv)}  if $u,\ v\in\mathfrak W$ and $u\overset{\mathfrak K}\approx v$, then $K\subset\mathbb N$ is $(u,\mathfrak K)$-small iff it is
$(v,\mathfrak K)$-small. \pf $|a_u(K,n)-a_v(K,n)|\le\sum_{k=0}^n|u_k-v_k|=o(a_u(n))$ along $\mathfrak K$.
\sms {\rm (v)}  if $u\in\mathfrak W$ is $\mathfrak K$-smooth, then $K\subset\mathbb N$ is $(u,\mathfrak K)$-small iff $K+1$ is
$(u,\mathfrak K)$-small.
\sms {\rm (vi)}  if $u\in\mathfrak W$ is $\mathfrak K$-smooth, $p\ge 1$, and $u^{(p)}_n:=u_{pn}$, then $a_{u^{(p)}}(n)\sim\tfrac1pa_u(pn)$ along $\tfrac1p\mathfrak K$ and
$u^{(p)}$ is $\tfrac1p\mathfrak K$-smooth where $\tfrac1p\mathfrak K:=\{\lfl\tfrac{j}p\rfl:\ j\in\frak K\}$.
\pf This follows from
$$\sum_{k=1}^n|u^{(p)}_k-u^{(p)}_{k+1}|\le p \sum_{k=1}^{pn}|u_k-u_{k+1}|=o(a_u(pn))\ \&$$
$$|pa_{u^{(p)}}(n)-a_u(pn)|\le\sum_{r=1}^{p-1}\sum_{k=1}^n|u_{pk}-u_{pk+r}|\le p^2\sum_{k=1}^{pn}|u_k-u_{k+1}|=o(a_u((p+1)n)).$$

\

\proclaim{Proposition 3.1}
\

 Suppose that $u\in\mathfrak W$ and $u_n\asymp \frac{ a_u(n)}n$ and  $u_n\asymp v_n$ where $v\in\mathfrak W$ and $v_n\downarrow$, then a set is  $u$-small  iff it has zero density.\endproclaim
\pf  By remark 3.1(ii), there is no loss of generality in assuming $u_n\downarrow$.

\pf of  $\Rightarrow$:
\

Suppose that $u_n\ge\eta\tfrac{ a_u(n)}n$. Since $u_n\downarrow$,
\begin{align*} a_u(K,n)=\sum_{k\in[1,n]\cap K}u_k \ge u_n|K\cap[1,n]|\ge \eta\frac{ a_u(n)}n|K\cap[1,n]|\end{align*}
 whence
 $$\frac{|K\cap[1,n]|}n\ \le\ \frac{ a_u(K,n)}{\eta  a_u(n)}.\ \ \CheckedBox$$

 \pf of $\Leftarrow$:
We show first that
$$\varlimsup_{n\to\infty,\ \e\to 0}\frac{ a_u(\e n)}{ a_u(n)}=0.$$
To do this, it suffices to show that
$$ a_u(2n)\gtrsim\ (1+\eta\log 2) a_u(n)$$
where  $u_n\ge\eta\tfrac{ a_u(n)}n$.
\

Indeed,
$$ a_u(2n)- a_u(n)=\sum_{k=n+1}^{2n}u_k\ge\eta\sum_{k=n+1}^{2n}\tfrac{ a_u(k)}k\ge\eta a_u(n)\sum_{k=n+1}^{2n}\tfrac1k\gtrsim\ \eta\log 2\cdot  a_u(n).\ \ \CheckedBox$$
Next, since $u_n\downarrow,\ u_n\le M\tfrac{ a_u(n)}n$, and
\begin{align*} a_u(K,n)=\sum_{k\in[1,n]\cap K}u_k & \le  a_u(n\e)+\sum_{k\in[\e n,n]\cap K}u_k\\ &\le  a_u(n\e)+u_{n\e}|[1,n]\cap K|\\ &\le
 a_u(n\e)+\frac{ a_u(n\e)}{n\e}|[1,n]\cap K|\end{align*}
whence
$$\frac{ a_u(K,n)}{ a_u(n)}\le \frac{ a_u(n\e)}{ a_u(n)}+\frac{1}{\e}\cdot\frac{|K\cap[1,n]|}n.\ \ \CheckedBox$$
\subsection*{Remark 3.2}
\

\sms (i) In case $u\in\mathfrak W$ and $u_n$ is {\it regularly varying} with {\it index} $s\in (-1,0)$
(i.e.\  $\tfrac{u_{\lfl\l n\rfl}}{u_n}\underset{n\to\infty}\lra\l^s$) then (see e.g. [BGT]) $u_n\sim \tfrac{(1+s) a_u(n)}n$ and $\exists\ v\in\mathfrak W$\f$\st$ $v_n\downarrow,\
v_n\sim u_n$. Thus, proposition 3.1  applies.

\sms (ii) The conclusion of proposition 3.1 fails for  $u_n=\tfrac1{n+1}$. The set
$K:=\bigcupdot_{k=1}^\infty [2^{k^2},k2^{k^2}]\cap\mathbb N$ is $u$-small, but $\varlimsup_{n\to\infty}\frac{|K\cap [1,n]|}n=1$.

\subsection*{$(u,\mathfrak K)$-density and $(u,\mathfrak K)$-strong Cesaro convergence}
\

Let $\mathfrak K\subset\mathbb N$ be a subsequence and let $u\in\mathfrak W$.
\

We'll say that a sequence $s_n$:
 \sbul {\it converges in $(u,\mathfrak K)$-density to} $L\in\mathbb R$ ($s_n\overset{(u,\mathfrak K)-\text{\tt\tiny d.}}{\underset{n\to\infty}\lra}\ L$) if
 $\exists\ K\subset\mathbb N$ $(u,\mathfrak K)$-small $\st$
$$s_n{\underset{n\to\infty,\ n\notin K}\lra}\ L;$$
and that
\sbul {\it $s_n$  converges   $(u,\mathfrak K)$-strongly Cesaro to} $L\in\mathbb R$ ($s_n\overset{(u,\mathfrak K)-\text{\tt\tiny s.C.}}{\underset{n\to\infty}\lra}\ L$) if

$$\tfrac1{a_u(n)}\sum_{k=0}^nu_k|s_k-L|{\underset{n\to\infty,\ n\in\frak K}\lra}\ 0.$$

\

\

\subsection*{Remark 3.3}
\

Let $\mathfrak K\subset\mathbb N$ be a subsequence and let $u\in\mathfrak W,\ s=(s_1,s_2,\dots)\in\mathbb R^\mathbb N$ and $L\in\mathbb R$.
\

\begin{align*}&\tag{i} s_n\overset{(u,\mathfrak K)-\text{\tt\tiny d.}}{\underset{n\to\infty}\lra}\ L\ \ \text{\tt if and only if}\\ &
K_\e:=\{n\in\mathbb N:\ |s_n-L|>\e\}\ \ \ \text{\tt is $(u,\mathfrak K)$-small}\ \forall\ \e>0.\end{align*}
\pf:\ \ Evidently $s_n\overset{(u,\mathfrak K)-\text{\tt\tiny d.}}{\underset{n\to\infty}\lra}\ L$ $\Rightarrow$ $K_\e$ is $(u,\mathfrak K)$-small $\forall\ \e>0$. To see  the reverse implication, if  $K_\e$ is $(u,\mathfrak K)$-small $\forall\ \e>0$, then by remark 3.1(iii) (on p. \pageref{rem3.1.3}),
 $\exists\ N_1<N_2<\dots$ so that $K:=\bigcup_{j=1}^\infty K_{1/j}\cap [N_j+1,N_{j+1}]$ is a $(u,\mathfrak K)$-small set and
 and $s_n\underset{n\to\infty,\ n\notin K}\lra 0$.
\sms (ii) if $v\in\mathfrak W,\ u\overset{\mathfrak K}\approx v$, then
\begin{align*}s_n\overset{(u,\mathfrak K)-\text{\tt\tiny d.}}{\underset{n\to\infty}\lra}\ L\ \ \text{\tt if and only if}\ \
s_n\overset{(v,\mathfrak K)-\text{\tt\tiny d.}}{\underset{n\to\infty}\lra}\ L.\end{align*}
\begin{align*}\tag{iii}  s_n{\overset{(u,\mathfrak K)-\text{\tt\tiny s.C.}}{\underset{n\to\infty}\lra}}\ L\ \Rightarrow\ \ s_n\overset{(u,\mathfrak K)-\text{\tt\tiny d.}}{\underset{n\to\infty}\lra}\ L.\end{align*}
\pf
\

 Suppose that $s_n\ge 0,\ L=0$ and that
$\tfrac1{ a_u(n)}\sum_{k=1}^n u_ks_k\underset{n\to\infty,\ n\in\mathfrak K}\lra\ 0$. By the Chebyshev-Markov inequality,
$$\frac{ a_u(K_\e,n)}{ a_u(n)}=\frac1{ a_u(n)}\sum_{k=1,\ k\in K_\e}^n u_k\le\frac1{\e  a_u(n)}\sum_{k=1}^n u_ks_k\underset{n\to\infty,\ n\in\mathfrak K}\lra\ 0.\ \ \CheckedBox$$
\

We call a sequence $x=(x_1,x_2,\dots)\in \mathbb R^\mathbb N$ {\it  one-sidedly bounded } if it is either bounded above, or below (or both).
\

\proclaim{Proposition 3.2}\ \  \ \ Let $\mathfrak K\subset\mathbb N$ be a subsequence and let $u\in\mathfrak W$.
Suppose that $x=(x_1,x_2,\dots)\in \mathbb R^\mathbb N$ is one-sidedly bounded, and $L\in\mathbb R$, then
\begin{align*}\tag{\ddag}x_n\overset{(u,\mathfrak K)-\text{\tt\tiny s.C.}}{\underset{n\to\infty}\lra}\ L\end{align*}
if and only if
\begin{align*}\tag{\dag}\frac1{ a_u(n)}\sum_{k=0}^nu_kx_k\underset{n\to\infty,\ n\in\mathfrak K}\lra\ L\ \&\ \ x_n\overset{(u,\mathfrak K)-\text{\tt\tiny d.}}{\underset{n\to\infty}\lra}\ L.\end{align*}

\endproclaim
\pf{\it of } (\dag) $\Rightarrow$ (\ddag):
\

We assume (without loss in generality) that $x_n\ge 0\ \forall\ n\ge 1$ and $L\ge 0$. Fix $\e>0$ and set
$$K_{+,\e}:=\{n\ge 1,\ x_n>L+\e\},\ K_{-,\e}:=\{n\ge 1,\ x_n<L-\e\},\ K_\e:=K_{+,\e}\cup K_{-,\e}.$$ By our assumptions $\exists\ N_\e\ge 1\ \st$
$$a_u(K_\e,[1,n])<\e a_u(n)\ \ \&\ \ \ L-\e<\frac1{ a_u(n)}\sum_{k=0}^nu_kx_k<L+\e\ \forall\ n>N_e,\ n\in\mathfrak K.$$
For large enough $n>N_\e,\ n\in\mathfrak K$,
\begin{align*}\sum_{k=1}^n&u_k|x_k-L| =(\sum_{k\in K_\e^c\cap [1,n]}+\sum_{k\in K_{+,\e}\cap [1,n]}+\sum_{k\in K_{-,\e}\cap [1,n]})u_k|x_k-L|\\ &\le
\e \sum_{k\in K_\e^c\cap [1,n]}u_k+\sum_{k\in K_{+,\e}\cap [1,n]}u_k(x_k-L)+\sum_{k\in K_{-,\e}\cap [1,n]}u_k(L-x_k)\\ &\le \e a_u(n)+La_u(K_\e,n)+\sum_{k\in K_{\e}\cap [1,n]}u_kx_k\\ &<
(1+L)\e a_u(n)+\sum_{k\in K_{\e}\cap [1,n]}u_kx_k
\end{align*}
Now (for large enough $n>N_\e,\ n\in\mathfrak K$),
\begin{align*}\sum_{k\in K_{\e}\cap [1,n]}u_kx_k& =\sum_{k\in [1,n]}u_kx_k-\sum_{k\in K_{\e}^c\cap [1,n]}u_kx_k\\ &<(L+\e)a_u(n)-(L-\e)a_u(K_\e^c,n)\\ &=
(L+\e)a_u(K_\e,n)+2\e a_u(K_\e^c,n)\\ &<\e(L+\e)a_u(n)+2\e a_u(n)\\ &=
(L+2+\e)\e a_u(n).\end{align*}
Reassembling,
$$\sum_{k=1}^nu_k|x_k-L|<(2L+3+\e)\e a_u(n).\  \ \CheckedBox$$
\

A version of the  following proposition is implicit in [GL]:
\proclaim{Proposition 3.3}
\

Let $\mathfrak K\subset\mathbb N$ be a subsequence and let $u\in\mathfrak W$.  Suppose that $x=(x_1,x_2,\dots)\in \mathbb R^\mathbb N$  and that $L\in\mathbb R$.

\

If $\frac1{ a_u(n)}\sum_{k=0}^nu_kx_k\underset{n\to\infty,\ n\in\mathfrak K}\lra\ L$ and
\par  either {\rm (i)} $x=(x_1,x_2,\dots)$  is bounded below and $\exists$ a $(u,\mathfrak K)$-small set $K_0\subset\mathbb N\ \st$
$\stackunder{n\to\infty,\ n\notin K_0}\varliminf\, x_n\ \ge\ \ L$;
\

or {\rm (ii)} $x=(x_1,x_2,\dots)$  is bounded above and  $\exists$ a $(u,\mathfrak K)$-small set $K_0\subset\mathbb N\ \st$
$\stackunder{n\to\infty,\ n\notin K_0}\varlimsup\, x_n\ \le\ \ L$;
\
then
$x_n\overset{(u,\mathfrak K)-\text{\tt\tiny s.C.}}{\underset{n\to\infty}\lra}\ L.$\endproclaim
\pf
\

By proposition 3.2, it suffices to prove that $x_n\overset{(u,\mathfrak K)-\text{\tt\tiny d.}}{\underset{n\to\infty}\lra}\ L$. By symmetry it suffices to prove the proposition under assumption (i). By possibly translating $x$ with a constant sequence, we reduce to the case $x_n\ge 0\ \forall\ n\ge 1\ \&\ L\ge 0$.
\

 For $\e>0$, set $K_\e:=\{n\notin K_0:\ x_n>L+\e\}$. It suffices to prove  that $K_\e$ is $(u,\mathfrak K)$-small $\forall\ 0<\e<\tfrac12$.

To see this, fix $0<\e<\tfrac12$ let $N_\e$ be so that
$$x_n>L-\e^2\ \forall\ n\ge N_\e,\ n\notin K_0,$$
then for large $n\gg N_\e,\ n\in\mathfrak K,$
\begin{align*}(L&+\e^2) a_u(n) >\sum_{k=N_\e}^nx_ku_k\\ &\ge \sum_{k\in K_\e\cap [N_\e,n]}x_ku_k+\sum_{k\in K_\e^c\cap K_0^c\cap[N_\e,n]}x_ku_k\\
&>(L+\e) a_u(K_\e,n)+(L-\e^2) a_u(K_\e^c\cap K_0^c,n)-(2L+\e-\e^2)a_u(N_\e)\\ &
>(L+\e) a_u(K_\e,n)+(L-\e^2) a_u(K_\e^c,n)- a_u(K_0,n)-(2L+1)a_u(N_\e)\\ &=(L+\e) a_u(K_\e,n)+(L-\e^2) a_u(K_\e^c,n)-\mathcal E_n\end{align*}
where
$\mathcal E_n:=a_u(K_0,n)+(2L+1)a_u(N_\e)$.
\

Writing $$(L+\e^2) a_u(n)=(L+\e^2) a_u(K_\e,n)+(L+\e^2) a_u(K_\e^c,n)$$ we see that
$$(\e-\e^2) a_u(K_\e,n)\le 2\e^2 a_u(K_\e^c,n)+\mathcal E_n.$$
For large $n\in\mathfrak K$,\  $\mathcal E_n<\e^2a_u(n)$ whence
$$ a_u(K_\e,n)\le \frac{3\e}{1-\e}\cdot  a_u(n)<6\e  a_u(n).\ \ \CheckedBox$$

\proclaim{Corollary 3.4}
\

Let $\mathfrak K\subset\mathbb N$ be a subsequence and let $u\in\mathfrak W$.  If $\tfrac{u_{n+1}}{u_n}\overset{(u,\mathfrak K)-\text{\tt\tiny d.}}{\underset{n\to\infty}\lra}\ 1$, then
$u$ is $\mathfrak K$-smooth.\endproclaim
\pf\ \ If $u_n>0\ \forall\ n\ge 0$, this follows from proposition 3.2 since
$\sum_{k=0}^n u_k|\tfrac{u_{k+1}}{u_k}-1|=\sum_{k=0}^n |u_{k+1}-u_k|.$
\

If this is not the case, define $v\in\mathfrak W$ by
$$v_n=\begin{cases} & u_n\ \ \ \ \ \ \ u_n>0,\\ & \frac1{2^n}\ \ \ \ \ \ \ u_n=0.\end{cases}$$
Evidently $a_u(n)\le a_v(n)\le a_u(n)+2$ so $a(n):=a_u(n)\sim a_v(n)$.
\

Moreover, if  $K\subset \mathbb N$, then
$a_u(K,n)\le a_v(K,n)\le a_u(K,n)+2$ and so $K$ is  $(u,\mathfrak K)$-small iff it is $(v,\mathfrak K)$-small.

Next, if $\tfrac12<\tfrac{u_{n+1}}{u_n}<2$, then $v_n=u_n$ and $v_{n+1}=u_{n+1}$,
$\tfrac{v_{n+1}}{v_n}=\tfrac{u_{n+1}}{u_n}$. Thus $\tfrac{v_{n+1}}{v_n}\overset{v-\text{\tt\tiny d.}}{\underset{n\to\infty}\lra}\ 1$ and by
 proposition 3.2 (as above),
 $$\tfrac1{a(n)}\sum_{k=0}^n |v_{k+1}-v_k|\underset{n\to\infty,\ n\in\mathfrak K}\lra 0.$$
 Finally,
 \begin{align*}\tfrac1{a(n)}\sum_{k=0}^n |u_{k+1}-u_k|&\le \tfrac1{a(n)}\sum_{k=0}^n |v_{k+1}-v_k|+\frac2{a(n)}\\ &\underset{n\to\infty,\ n\in\mathfrak K}\lra 0.\ \ \CheckedBox\end{align*}

\section*{\S4 Proofs of proposition 0 and theorem A}
\demo{Proof of proposition 0}
\

Fix $F\in\mathcal F$ satisfying ($\bigstar_\mathfrak K$). Evidently, for $A,\ B\in\B\cap F$,

$$\frac1{a_n(F)}\sum_{k=0}^{n-1}m(A\cap T^{-k}B)\underset{n\to\infty,\ n\in\mathfrak K}\lra\ m(A)m(B)\ \forall\ A,\ B\in\B\cap F.$$
 This shows that $F\in R_\frak K(T)$ and that $T$ is weakly rationally ergodic along $\mathfrak K$; proving (i).

\

To prove (ii), let $F\in\mathcal F$ satisfy ($\bigstar_\mathfrak K$). It suffices to show  that
 \begin{align*}\tag{\dsbiological}m(B\cap T^{-n}C)\ \overset{\mathfrak K}\approx\ m(B)m(C)u_n(F)\ \forall\ B,\ C\in R_{\mathfrak K}(T).\end{align*}
 \pf {\it of} (\dsbiological):
 \

 Fix $B,\ C\in R_{\mathfrak K}(T)$,  then $G:=B\cup C\in R_{\mathfrak K}(T)$ and we claim:
 \Par1 $\exists\ K\subset\mathbb N$ $(u(F),\mathfrak K)$-small $\st$
$$\varliminf_{n\to\infty,\ n\notin K}\frac{m(B\cap T^{-n}C)}{u_n(F)}\ge m(B)m(C).$$
\

\pf {\it of} \P1:
  \

  Let $\e>0$, then $\exists\ B_0,\dots ,B_N,\ C_0,\dots ,C_N\in\B\cap F$ so that $$B':=\bigcupdot_{k=0}^NT^{-k}B_k\subset B,\  C':=\bigcupdot_{k=0}^NT^{-k}C_k\subset C,\
m(B\setminus B')<\e,\ \ m(C\setminus C')<\e.$$
Using ($\bigstar_{\mathfrak K}$) (as on page  \pageref{bigstarK})
\begin{align*}\frac{m(B\cap T^{-n}C)}{u_n(F)}&\ge \frac{m(B'\cap T^{-n}C')}{u_n(F)}\\ &=\sum_{k,\ell=0}^N\frac{m(T^{-k}B_k\cap T^{-n-\ell}C_\ell)}{u_n(F)}\\ &\overset{(u(F),\mathfrak K)-\text{\tt\tiny d.}}{\underset{n\to\infty}\lra}\sum_{k,\ell=0}^Nm(T^{-k}B_k)m( T^{-\ell}C_\ell)\\ &=m(B')m(C')>(m(B)-\e)(m(C)-\e).\end{align*}
Choose $\e_n\downarrow\ 0$. By the above $\exists\ K_1\subset K_2\subset\dots\subset\Bbb N$, each $K_\nu$ being $(u(F),\mathfrak K)$-small, $\st$
$$\varliminf_{n\to\infty,\ n\notin K_\nu}\frac{m(B\cap T^{-n}C)}{u_n(F)}\ge m(B)m(C)-\e_\nu\ \ \forall\ \nu\ge 1.$$
By remark 3.1(iii) (on p. \pageref{rem3.1.3}) $\exists\ K\subset\Bbb N$ realizing  \P1.
\

By weak rational ergodicity along $\mathfrak K$,
$$\frac1{a_n(F)}\sum_{k=0}^nm(B\cap T^{-k}C){\underset{n\to\infty,\ n\in\mathfrak K}\lra} m(B)m(C)$$
so by proposition 3.3(i),
$$\frac{m(B\cap T^{-n}C)}{u_n(F)}\ \overset{(u(F),\mathfrak K)-\text{\tt\tiny s.C.}}{\underset{n\to\infty}\lra} m(B)m(C)\ \ \forall\ B,\ C\in \B\cap A;$$ equivalently
$m(B\cap T^{-n}C)\overset{\mathfrak K}\approx\ m(B)m(C)u_n(F).\ \ \CheckedBox\text{(\dsbiological)}$

 This proves (ii).
\

To see (iii) ($\mathfrak K$-smoothness of $u(F)$ for $F\in R_{\mathfrak K}(T))$),  take $B=F,\ C=T^{-1}F$ in (\dsbiological).
\

To prove (iv), fix $p\in\Bbb N$. To see rational weak mixing of $T^p$ along $\tfrac1p\frak K$, let $A\in R_\frak K(T)$.
By  remark 3.1(vi),
$a_n^{(T^p)}(A)\sim\tfrac1pa_{np}^{(T)}(A)$ along $\tfrac1p\frak K$. It also follows that for $B,C\in\B\cap A$,
\begin{align*}\sum_{k=0}^n|m(B\cap T^{-kp}C)&-m(B)m(C)u^{(p)}_k|=\sum_{k=0}^n|m(B\cap T^{-kp}C)-m(B)m(C)u_{pk}|\\ &
\le\sum_{k=0}^{pn}|m(B\cap T^{-k}C)-m(B)m(C)u_{k}|\\ &=o(a_{pn}(A))\ \ \ \text{\rm along}\ \ \tfrac1p\frak K.
\end{align*} This shows that $R_\frak K(T)\subset R_{\tfrac1p\frak K}(T^p)$.
\

The other inclusion follows from results in [FL]. The proof of theorem 3.3 there shows that
\begin{align*}A\in R_{\frak K}(T)\ \Leftrightarrow\ \ \{\tfrac{S_n^{(T)}(1_A)}{a_n^{(T)}(A)}:\ n\in\frak K\}\ \ \text{\tt is uniformly integrable on}\ A.\end{align*}

Now
$$S_{pn}^{(T)}(1_A)=\sum_{k=0}^{pn-1}1_A\circ T^{k}=\sum_{\nu=0}^{p-1}\sum_{k=0}^{n-1}1_A\circ T^{kp+\nu}=\sum_{\nu=0}^{p-1}S_{n}^{(T^p)}(1_A)\circ T^{\nu}$$
whence
\begin{align*}A\in R_{\frac1p\frak K}(T^p)&\ \Leftrightarrow\ \ \{\tfrac{S_n^{(T^p)}(1_A)}{a_n^{(T^p)}(A)}:\ n\in\frac1p\frak K\}\ \ \text{\tt is uniformly integrable on}\ A
\\ &\Lra\ \ \ \{\tfrac{S_{pn}n^{(T)}(1_A)}{a_{pn}^{(T)}(A)}:\ n\in\frac1p\frak K\}\ \ \text{\tt is uniformly integrable on}\ A
\\ &\Lra\ \ \ \{\tfrac{S_n^{(T)}(1_A)}{a_n^{(T)}(A)}:\ n\in\frak K\}\ \ \text{\tt is uniformly integrable on}\ A\\ &\Lra\ A\in R_\frak K(T).\ \ \ \CheckedBox\end{align*}
\ \ \ \CheckedBox

\demo{Proof of theorem A}\ \ \ This follows from proposition 0(ii) via proposition 3.1.\ \ \CheckedBox

\section*{\S5 Proof of Lemmas B and C, and corollary E}
\demo{Proof of Lemma B}

Let $\mathcal U$ be the collection of finite unions  sets in $\mathcal C$. It follows from assumptions (ii) (a) and (b) that

\begin{align*}\tag{\Football}\frac{m(A\cap T^{-k}B)}{u_k}\overset{(u,\mathfrak K)-\text{\tt\tiny d.}}{\underset{k\to\infty}\lra }  m(A)m(B)\ \ \forall\ A,\
 B\in\mathcal U.
\end{align*}

Let
$$\mathcal K:=\{K\subset\Om:\ \ K\ \text{\tt compact}\}.$$

We claim first that
\begin{align*}\tag{\Yinyang}\frac{m(A\cap T^{-k}B)}{u_k}\overset{(u,\mathfrak K)-\text{\tt\tiny d.}}{\underset{k\to\infty}\lra }  m(A)m(B)\ \ \forall\ A,\ B\in\mathcal K,\ \ A,\ B\subset\Om.
\end{align*}
\pf We show first that $\forall\ A,\ B\in\mathcal K,\ \ \exists\ K_1\subset\mathbb N$ $(u,\mathfrak K)$-large,  so that
\begin{align*}\tag{1}\varlimsup_{k\to\infty,\ k\in K_1}\frac{m(A\cap T^{-k}B)}{u_k}\le m(A)m(B).
\end{align*}\ \ \ To see this, we  show first that $\forall\ \e>0,\ \exists\ U,\ V\in\mathcal U$ so that
\begin{align*}\tag{2}A\subset U,\ B\subset V,\ m(U\setminus A),\ m(V\setminus B)<\e.\end{align*}
Given $\e>0$, the Borel property of the measure $m$ ensures  open sets $U,\ V$ satisfying (2). Each of these is a countable union of members of $\mathcal C$. By compactness of $A,\ B$ we can reduce to finite unions  $U,\ V\in\mathcal U$.
\

By (\Football) $\exists\ K_\e\subset\mathbb N$ $(u,\mathfrak K)$-small,  so that
\begin{align*} \frac{m(U\cap T^{-k}V)}{u_k}\underset{k\to\infty,\ k\notin K_\e}\lra m(U)m(V)
\end{align*}
and
\begin{align*}\varlimsup_{k\to\infty,\ k\in K_\e}\frac{m(A\cap T^{-k}B)}{u_k}< (m(A)+\e)(m(B)+\e).
\end{align*}
Fix $N_\nu\ \ (\nu\ge 1)$ so that $ a_u(K_{\frac1\nu},n)<\tfrac1\nu a_n(\Om)\ \forall\ n\in\mathfrak K,\ n\ge\nu$ and
$$\frac{m(A\cap T^{-k}B)}{u_k}< (m(A)+\frac1\nu)(m(B)+\frac1\nu)\ \forall\ k\in K_{\frac1\nu}^c\cap [N_\nu,\infty).$$
The set $$K_1:=\bigcupdot_{\nu\ge 1}K_{\frac1\nu}\cap [N_\nu,N_{\nu+1})$$
is as required for (1).
\

 Now fix $\forall\ A,\ B\in\mathcal K,\ \ A,\ B\subset\Om$. Since   $\Om\in R_{\mathfrak K}(T)$,

$$\sum_{k=0}^{n-1}m(A\cap T^{-k}B)\sim m(A)m(B)a_n(\Om)\ \ \ \text{\tt as}\ n\to\infty,\ n\in\mathfrak K,$$
and the claim follows from (1) and proposition 3.3(ii).\ \ \ \CheckedBox(\Yinyang)

\

To complete the proof of that $\Om$ satisfies ($\bigstar_\mathfrak K$), let $A,\ B\in\B\cap \Om$, then $\exists$ $E_N,\ F_N\in\mathcal K\ \st\ \mod m$:
$$E_N\uparrow A\ \ \&\ \ \ F_N\uparrow B$$
whence by  (\Yingyang), $\exists\ K_N$ $(u,\mathfrak K)$-small so that $\forall\ N\ge 1$,
$$\varliminf_{n\to\infty,\ n\notin K_N}\tfrac{m(A\cap T^{-n}B)}{u_n}\ge m(E_N)m(F_N).$$
As above, $\exists\ K\subset\ \mathbb N$ $(u,\mathfrak K)$-small $\st$
$$\varliminf_{n\to\infty,\ n\notin K}\tfrac{m(A\cap T^{-n}B)}{u_n}\ge m(A)m(B)$$and ($\bigstar_\frak K$) follows from proposition 3.3(i).\ \ \CheckedBox
\demo{Proof of Lemma C}

  By standardness, up to isomorphism, $X=\a^\mathbb Z$,  $T:X\to X$ is the shift and the collection $\mathcal C_\a$ of $(\a,T)$-cylinder sets forms a
base of clopen sets for the Polish topology on $X$. Thus $\mathcal C:=\mathcal C_\a\cap\Om$ satisfies assumptions (ii) of lemma B and lemma C follows.\ \ \ \CheckedBox

\demo{Proof  of Corollary E}

 Let $(X,\B,m,T)$ be  rationally weakly mixing  along $\mathfrak K$, and let $\pi:(X',\B',m',T')\to (X,\B,m,T)$ be its natural extension, that is:
$$T'\ \text{\tt\small invertible,}\ \pi\circ T'=T\circ\pi,\  \pi^{-1}\B\subset\B',\ m'\circ\pi^{-1}=m\ \&\  \bigvee_{n\ge 0}T^{\prime n}\pi^{-1}\B=\B'.$$
It follows  from   uniform integrability considerations (as in theorem 3.3 of [FL]) that $T'$ is weakly rationally ergodic  along $\mathfrak K$ with \f$R_{\mathfrak K}(T')\supseteq\mathcal{HR}(\pi^{-1}R_{\mathfrak K}(T))$.
\

To see that $T'$ is rationally weakly mixing along $\mathfrak K$, fix a countable, one-sided $T$-generator $\a\subset R_{\mathfrak K}(T)$, then $\a':=\pi^{-1}\a\subset R_{\mathfrak K}(T')$ is a countable, two-sided  $T'$-generator.
\

Fix $\Om\in \a'$ and let $u:=u(\Om)$. By rational weak mixing of $T$ along $\mathfrak K$,
\begin{align*}\frac{m'(A\cap T^{\prime -n}B)}{u_n}\ \ \overset{(u,\mathfrak K)-\text{\tt\tiny d.}}{\underset{n\to\infty}\lra}\ \ \
 m'(A)m'(B)\ \ \forall\ A,\ B\in\mathcal C_{\a'}=\pi^{-1}\mathcal C_{\a}\end{align*} whence by lemma C,
 $(X',\B',m',T')$ is rationally weakly mixing along $\mathfrak K$.\ \  \CheckedBox

\section*{\S6 Mean ergodic theorem for weighted averages}
\

Let $\mathfrak K\subset\mathbb N$ be a subsequence.  Call a weight $u\in\mathfrak W$ (good for the) {\it mean ergodic theorem along $\mathfrak K$} (abbr. {\tt MET}$_\mathfrak K$) if for any ergodic,  probability preserving transformation $(\Om,\mathcal A,P,S)$, we have that
\begin{align*}\tag{\tt MET$_\mathfrak K$}\tfrac1{a_u(n)}\sum_{k=0}^{n-1}u_kf\circ S^k\overset{L^2(P)}{\underset{n\to\infty,\ n\in\mathfrak K}\lra}\ E(f)\ \ \forall\ f\in L^2(P).\end{align*}
We let {\tt MET}:={\tt MET}$_\mathbb N$.

\

Let $\mathfrak K\subset\mathbb N$ be a subsequence.

Using the spectral theorem for unitary operators, it can be shown (see [Kre]) that the following conditions are equivalent for $u\in\mathfrak W$:
\bul $u$ is {\tt MET}$_\mathfrak K$;\bul $\frac1{a_u(n)}\sum_{k=0}^nu_kz^k\underset{n\to\infty,\ n\in\mathfrak K}\lra 0$  for $z\in\mathbb C,\ |z|=1,\ z\ne 1$;
\bul $u$ is (good for the) {\it  weak ergodic theorem along $\mathfrak K$} \ \ \ (abbr. {\tt WET}$_\mathfrak K$ )  in the sense that for any ergodic,  probability preserving transformation $(\Om,\mathcal A,P,S)$,
 {\small\begin{align*}\tag{\tt WET$_\mathfrak K$}\tfrac1{a_u(n)}\sum_{k=0}^{n-1}u_kf\circ S^k\underset{n\to\infty,\ n\in\mathfrak K}\lra\ E(f)\ \text{\tt\scriptsize  weakly in}\ L^2(P)\ \forall\  f\in  L^2(P).\end{align*}}

The  recurrent, renewal sequences form an important subclass of weights. A weight $u\in\mathfrak W$ is a {\it recurrent renewal sequence} if $u_0=1$ and $\exists\ f\in\mathcal P(\mathbb N)$, called the  associated {\it lifetime distribution}  satisfying the {\it renewal equation} $$u_n=\sum_{k=1}^nf_ku_{n-k}\ \ (n\ge 1).$$
\

The renewal sequence is $u$ is called {\it aperiodic} if $\<\{n\in\mathbb N:\ u_n>0\}\>=\mathbb Z$.
\

It follows from the renewal equation that any aperiodic, recurrent renewal sequence satisfies
$|\sum_{k=0}^nu_kz^k|<\infty$  for $z\in\mathbb C,\ |z|=1,\ z\ne 1$ and hence is {\tt MET}. Proposition 6.2 (below) generalizes this.
\

Let $\mathfrak K\subset\mathbb N$ be a subsequence. Any $\mathfrak K$-smooth weight $u\in\mathfrak W$ is {\tt MET}$_\mathfrak K$ (see [HP], [Kre] and references therein).
A weight $u\in\mathfrak W$  which is {\tt MET} and not $\mathfrak K$-smooth for  any subsequence $\mathfrak K\subset\mathbb N$ is exhibited in [HP].
 \

We'll need
\proclaim{Lemma 6.1}
\

Let $\mathfrak K\subset\mathbb N$ be a subsequence and suppose that $u\in\mathfrak W$ is {\tt MET}$_\mathfrak K$, and that $(\Om,\mathcal A,P,S)$ is a weakly mixing probability preserving transformation, then
$$P(A\cap S^{-n}B)\overset{(u,\mathfrak K)-\text{\tt\tiny d.}}{\underset{n\to\infty}\lra}\ P(A)P(B)\ \ \forall\ A,\ B\in\mathcal A.$$\endproclaim\pf
\

 It follows from ({\tt WET}$_\mathfrak K$) for $S$ and $A,\ B\in\mathcal A$, that
\begin{align*}\tag{\dsmilitary}\tfrac1{a_u(n)}\sum_{k=0}^{n-1}u_kP(A\cap S^{-k}B){\underset{n\to\infty,\ n\in\mathfrak K}\lra}\ P(A)P(B),\end{align*}
and it follows from (\dsmilitary) for $S\x S$ (which is ergodic) and $A\x A,\ B\x B\in\mathcal A\otimes\mathcal A$ that
\begin{align*}\tag{\dsheraldical}\tfrac1{a_u(n)}\sum_{k=0}^{n-1}u_kP(&A\cap S^{-k}B)^2\\ &=
\tfrac1{a_u(n)}\sum_{k=0}^{n-1}u_kP\x P(A\x A\cap (S\x S)^{-k}B\x B)\\ &{\underset{n\to\infty,\ n\in\mathfrak K}\lra}\ P\x P(A\x A)\ \cdot\ P\x P(B\x B)\\ &= P(A)^2P(B)^2.\end{align*}
Using (\dsmilitary) and (\dsheraldical)
\begin{align*}&\tfrac1{a_u(n)}\sum_{k=0}^{n-1}u_k(P(A\cap S^{-k}B)-P(A)P(B))^2=\\ &
\tfrac1{a_u(n)}\sum_{k=0}^{n-1}u_k\(P(A\cap S^{-k}B)^2-2P(A)P(B)P(A\cap S^{-k}B)+ P(A)^2P(B)^2\)\\ &
{\underset{n\to\infty,\ n\in\mathfrak K}\lra}\ 0\end{align*}
whence $P(A\cap S^{-n}B)\overset{(u,\mathfrak K)-\text{\tt\tiny d.}}{\underset{n\to\infty}\lra}\ P(A)P(B)$.\ \ \CheckedBox

\proclaim{Proposition 6.2}
\

Let $\mathfrak K\subset\mathbb N$ be a subsequence.  Suppose that $(X,\B,m,T)$ is weakly rationally ergodic along $\mathfrak K$ and spectrally weakly mixing, then
$u(E,F)$ is {\tt MET}$_\mathfrak K$ $\forall\ E,\ F\in R_{\mathfrak K}(T)_+$.\endproclaim\pf
\

Let $(\Om,\mathcal A,P,S)$ be an ergodic,  probability preserving transformation.

\

It follows from the assumptions that
$T\x S$ is weakly rationally ergodic along $\mathfrak K$  and $R_{\mathfrak K}(T\x S)\supset R_{\mathfrak K}(T)\x\Om$.
\

It suffices to show that for  $E,\ F\in R_{\mathfrak K}(T)_+$,
\begin{align*}\tag{\tt WET$_\mathfrak K$}A_nf\underset{n\to\infty,\ n\in\mathfrak K}\lra\ E(f) \ \text{\tt weakly in}\ L^2(P)\ \forall\  f\in  L^2(P)\end{align*}
where $A_nf:= \frac1{a_n(E,F)}\sum_{k=0}^{n-1}u_k(E,F)f\circ S^k$.

Since $E\x\Om,\ F\x\Om\in R_{\mathfrak K}(T\x S)$,
\begin{align*}\frac1{a_n}\sum_{k=0}^{n-1}u_k(E,F)P(C\cap S^{-k}D)\underset{n\to\infty}\lra\ P(C)P(D)\ \ \forall\ \ C,D\in\B(\Om).\end{align*}
This shows ({\tt WET}$_\mathfrak K$) for indicators, whence for simple functions $f$. By the  triangle inequality
$\|A_nf\|_2\le \|f\|_2 \forall\  f\in  L^2(P)$ and ({\tt WET}$_\mathfrak K$) follows by approximation .\ \ \ \CheckedBox

\section*{\S7 Proof of theorem D}
We assume that $T$ is invertible. By Corollary E, this involves no loss in generality.
\demo{Proof  of theorem D(i):}
\

 We'll prove spectral weak mixing of $T$ by showing that $T\x S$ is weakly rationally ergodic along $\mathfrak K$ for any ergodic, probability preserving transformation
 $(\Om,\mathcal A,P,S)$. To this end, let $(\Om,\mathcal A,P,S)$ be an ergodic, probability preserving transformation.  We claim  first that

 \

\begin{align*}\tag{\symrook} \  \ \frac1{a_n(T)}&\sum_{k=0}^{n-1}m(A\cap T^{-k}B)P(C\cap S^{-k}D)\\ &\xrightarrow[n\to\infty]{ n\in\mathfrak K}\ m(A)m(B)P(C)P(D)\ \forall\ A,\ B\ \in R_{\mathfrak K}(T),\ C,D\in\mathcal A.\end{align*}
\pf {\it of} \ (\symrook):
\

Fix $A,\ B\ \in R_{\mathfrak K}(T)_+$ and set $v=u(A,B)$. By proposition 0(iii), $v$ is smooth whence {\tt MET}$_\frak K$; and (\symrook) follows from ({\tt WET}$_\frak K$) for $C,\ D\in\mathcal A$.\ \ \CheckedBox
\

\

Let $\mu:=m\x P,\ \C:=\B\otimes\mathcal A$ and $\tau:=T\x S$.
\

 We claim next that for $ F\in R_{\mathfrak K}(T),\ A\in\B\cap F,\  B\in\mathcal A,\ C\in\C\cap (F\x \Om)$,
\begin{align*}\tag{\symbishop} \tfrac1{a_n}\sum_{k=0}^{n-1}\mu(C\cap \tau^{-k}(A\x B))\underset{n\to\infty,\ n\in\mathfrak K}\lra\ \mu(C)\mu(A\x B).\end{align*}
\pf {\it of} \ (\symbishop):
 \

 By weak rational ergodicity  along $\mathfrak K$ the collection $\{\tfrac1{a_n}\sum_{k=0}^{n-1}1_A\circ T^k:\ n\in\mathfrak K\}$ is uniformly integrable on $F$.
It follows that the collection $\{\tfrac1{a_n}\sum_{k=0}^{n-1}1_{A\x B}\circ \tau^k:\ n\in\mathfrak K\}$ is uniformly integrable on $F\x\Om$.
\

Let $\Phi\in \{\tfrac1{a_n}\sum_{k=0}^{n-1}1_{A\x B}\circ \tau^k:\ n\in\mathfrak K\}'$ be a weak limit, then by (\symrook),
$$\int_{C\x D}\Phi d\mu=\mu(A\x B)\mu(C\x D)\ \forall\ C\in\B\cap F,\ D\in\mathcal A.$$
It follows that $\Phi\equiv\mu(A\x B)$ whence $\tfrac1{a_n}\sum_{k=0}^{n-1}1_{A\x B}\circ \tau^k\underset{n\to\infty,\ n\in\mathfrak K}\lra\ \mu(A\x B)$ weakly in $L^1(F\x\Om)$ and (\symbishop) follows.\ \ \CheckedBox
\

Finally we complete the proof of theorem D(i) by showing that
$$F\x\Om\in R_{\mathfrak K}(\tau);$$
namely,  for
$ F\in R_{\mathfrak K}(T),\ C, D\in\C\cap (F\x \Om)$,
\begin{align*}\tag{\symknight} \frac1{a_n(F)}\sum_{k=0}^{n-1}\mu(C\cap \tau^{-k}D)\underset{n\to\infty,\ n\in\mathfrak K}\lra\ \mu(C)\mu(D).\end{align*}

\

\pf{\it of} (\symknight):
\

Since $D\subset F\x\Om$, the collection $\{\tfrac1{a_n}\sum_{k=0}^{n-1}1_D\circ \tau^k:\ n\in\mathfrak K\}$ is uniformly integrable on $F\x\Om$. Let $\Psi\in \{\tfrac1{a_n}\sum_{k=0}^{n-1}1_{D}\circ \tau^k:\ n\in\mathfrak K\}'$ be a weak limit, then by (\symbishop) for $\tau^{-1}$,
$$\int_{A\x B}\Phi d\mu=\mu(A\x B)\mu(D)\ \forall\ A\in\B\cap F,\ B\in\mathcal A$$
 whence $\tfrac1{a_n}\sum_{k=0}^{n-1}1_{D}\circ \tau^k\underset{n\to\infty,\ n\in\mathfrak K}\lra\ \mu(D)$ weakly in $L^1(F\x\Om)$ and (\symknight) follows.\ \ \CheckedBox
\

\subsection*{Remark}\ \ Spectral  weak mixing alone does not imply subsequence rational weak mixing. See [ALV] for {\tt squashable}, spectrally weakly mixing, transformations. These are not even subsequence  weakly rationally ergodic.
 We do not know whether weak rational ergodicity and spectral weak mixing together imply subsequence rational weak mixing.
\
 \demo{Proof  of theorem D(ii):}
\

Fix a countable, $\B$-generating partition $\a\subset R_{\mathfrak K}(T)$.
By standardness, up to isomorphism, $X=\a^\mathbb Z$ and $T:X\to X$ is the shift.
The collection $\mathcal C_\a$ of $(\a,T)$-cylinder sets forms a base of clopen sets for the $T$-invariant, measurable, Polish topology on $X$.
\

Let $(\Om,\mathcal A,P,S)$ be a weakly mixing, probability preserving transformation.
 Fix a compact $S$-invariant, completely disconnected, measurable topology on $\Om$ generating $\mathcal A$.
\

We must show that the  measure preserving transformation
$$(Z,\C,\mu,\tau):=(X\x\Om,\B\otimes\mathcal A,m\x P,T\x S)$$ is rationally weakly mixing along $\mathfrak K$.
\

For this, it suffices to show that for  $F\in \a,\ m(F)>0$,\ \ $F\x\Om$ satisfies ($\bigstar_\mathfrak K$) $\wrt$ $\tau $.
\

By Lemma B, it suffices to establish
\begin{align*}\tag{\symking} \frac{m(A\cap T^{-n}B)P(C\cap S^{-n}D)}{u_n(F)}&\overset{(u(F),\mathfrak K)-\text{\tt\tiny d}}{\underset{n\to\infty}\lra}\ m(A)m(B)P(C)P(D)\\ &\ \forall\ A,\ B\in\B\cap F,\ C,D\in\mathcal A.\end{align*}
\

\pf {\it of} (\symking):
\

 By proposition 0(iii), $u(F)$ is  $\frak K$-smooth, whence {\tt MET}$_\mathfrak K$ and by lemma 6.1
$$P(C\cap S^{-n}D)\overset{(u(F),\mathfrak K)-\text{\tt\tiny d}}{\underset{n\to\infty}\lra}\ P(C)P(D)\ \ \forall\ \ C,D\in\mathcal A.$$
Thus, since $F$ satisfies ($\bigstar_\mathfrak K$), for $A,\ B\in\B\cap F$
$$\frac{m(A\cap T^{-n}B)}{u_n(F)}\overset{(u(F),\mathfrak K)-\text{\tt\tiny d}}{\underset{n\to\infty}\lra}\ m(A)m(B).$$
These two $(u(F),\mathfrak K)$-density convergences imply (\symking), and (via lemma B) theorem D(ii).\ \ \ \CheckedBox

\section*{\S8  Markov shift examples}
\

\par Let $S$ be a countable set (the {\it state space}) and let  $P:S\x S\to [0,1]$ be a {\tt stochastic matrix}  (the {\it transition matrix}) on $S$ ($\sum_{t\in S}p_{s,t}=1\
\forall s\in S$) with an invariant distribution $\pi:S\to \mathbb R_+$ (${\sum_{u\in S}}\pi_up_{u,t}=\pi_t$).

\par The stationary, {\it two-sided  Markov shift} of $(P,\pi)$
is the quadruple $$(S^{\mathbb Z},{\mathcal  B},m,T),$$ where  $T:S^{\mathbb Z}\to S^{\mathbb Z}$ is the
 shift,
$$\B:=\s(\{\text{\tt cylinders}\}),$$  a {\it cylinder } being a set of form

$$[s_1,\dots,s_n]_k:=\{x=(\dots,x_{-1},x_0,x_1,\dots)\in S^{\mathbb
Z}:x_{j+k}=s_{j}\ \forall\ 1\le j\le n\}$$ $(s_1,\dots,s_n\in S^n,\ k\in\mathbb Z,\ n\in\mathbb N)$;  and the measure $m$ is defined by
$$m([s_1,\dots,s_n]_k)=\pi_{s_1}p_{s_1,s_2}\cdots p_{s_{n-1},s_n}\ \forall
s_1,\dots,s_n\in S^n,\ n\in\mathbb N.$$

\

 The stationary Markov shift $(S^{\mathbb Z},{\mathcal  B},m,T),$  is a measure preserving transformation.
 \

 As shown in [HR],  $T$ is\sbul conservative iff $P$ is {\it recurrent}
 ($\sum_{n=0}^\infty p_{s,s}^{(n)}=\infty \ \forall\ s\in S$)
  \

  and in this case, $T$ is \sbul ergodic
 iff $P$ is  {\it irreducible} ($\forall s,t\in S,\ \ \exists n\in\mathbb N\ \ni p_{s,t}^{(n)}>0$).

 \

The (stationary) {\it one-sided, Markov shift} is  $(S^{\mathbb N},{\mathcal  B_+},m_+,\tau),$   where
$\tau:S^{\mathbb N}\to S^{\mathbb N}$ is the
 shift, $$\B_+:=\s(\{\text{\tt one-sided cylinders}\}),$$  a {\it one-sided cylinder } being a set of form
$$[s_1,\dots,s_n]:=\{x=(x_1,x_2,\dots)\in S^{\mathbb
N}:x_{j}=s_{j}\ \forall\ 1\le j\le n\}$$ $(s_1,\dots,s_n\in S^n,\ k\in\mathbb Z,\ n\in\mathbb N)$;  and the measure $m_+$ is defined by
$$m_+([s_1,\dots,s_n])=\pi_{s_1}p_{s_1,s_2}\cdots p_{s_{n-1},s_n}\ \forall
s_1,\dots,s_n\in S^n,\ n\in\mathbb N.$$

\

 As shown in [BF], if the stochastic matrix $P$ is irreducible, recurrent and {\it aperiodic}
 ($\gcd\,\{n\ge 1:p_{s,s,}^{(n)}>0\}=1$ for some and hence all $s\in S$), then $T$ is a conservative K-automorphism (natural extension of an exact endomorphism), whence (see [ALW]) spectrally weakly mixing.
 \

 As shown in [A1], a conservative, ergodic Markov shift $(S^{\mathbb N},{\mathcal  B},m,T)$ is rationally ergodic with $R_{\mathfrak K}(T)\supset\mathcal{HR}(\mathcal C_\a)$
 where $\a:=\{[s]_0:\ s\in S\}$;   with
 $a_n(T)=a_n(P)\sim \tfrac1{\pi_s}\sum_{k=0}^{n-1}p_{s,s}^{(k)}\ \ \ (\forall\ s\in S)$.

\proclaim{Theorem 8.1}
\

Let $\mathfrak K\subset\mathbb N$ be a subsequence. The Markov shift $(S^{\mathbb Z},{\mathcal  B},m,T)$   of the irreducible, recurrent, aperiodic transition  matrix
$P$ on state space $S$ is  rationally weakly mixing along $\mathfrak K$ iff
 $\exists\ s\in S$ with
$u([s]_0)$
  $\mathfrak K$-smooth.
\endproclaim
\pf By proposition 0(iii), if $T$ is  rationally weakly mixing along $\mathfrak K$, then  $u([s]_0)$ is $\mathfrak K$-smooth $\ \forall\ s\in S$.
 \

 To prove the other implication, we'll need the following lemma:

 \

\proclaim{Lemma 8.2}
\

Let $S$ be a countable set and let $P:S\x S\to [0,1]$ be an irreducible, recurrent, aperiodic stochastic matrix with the property that for some $s\in S$, $u=u([s]_0)$
is $\mathfrak K$-smooth, then

\begin{align*}\tag{\Bouquet}
\frac{p_{r,t}^{(n+\ell)}}{u_n}
\overset{(u,\mathfrak K)-\text{\tt\tiny s.C.}}{\underset{n\to\infty}\lra}\ \pi_t\ \ \forall\ r,\ t\in S,\ \ell\in\mathbb Z.\end{align*}\endproclaim
Lemma 8.2 is a $(u,\mathfrak K)$-density version of lemma 1 in [O].

\

\pf{\it of lemma 8.2}:
\

Recall from  [Ch] that the $P$-stationary distribution $\pi:S\to\mathbb R_+$ with $\pi_s=1$ is given by
$$\pi_t={\sum_{n=1}^\infty} \,_sp^{(n)}_{s,t}$$ where
$$\,_sp^{(1)}_{s,t}:=p_{s,t},\ \ \&\ \  \,_sp^{(n+1)}_{s,t}:=
{\sum_{r\in S\setminus\{s\}}} {_sp_{s,r}^{(n)}p_{r,t}}.$$

As shown in [Ch], $\forall\ r,\ t\in S$,
$$\frac1{a_u(n)}\sum_{k=0}^{n-1}p_{r,t}^{(k)}\underset{n\to\infty}\lra\ \pi_t.$$
In view of this, to show (\Bouquet), it suffices  by proposition 3.3(i)  to show that $\forall\ r,\ t\in S,\ \ell\in\mathbb Z$
$\exists\ K_{r,t,\ell}\subset\mathbb N$, $(u,\mathfrak K)$-small\ $\st$
\begin{align*}\tag{\dsagricultural} \stackunder{n\to\infty,\ n\notin K_{r,t,\ell}}\varliminf\, \frac{p_{r,t}^{(n+\ell)}}{u_n}\ \ge\ \ \pi_t.\end{align*}
Let $K_0\subset\mathbb N$ be $(u,\mathfrak K)$-small $\st$ $\frac{u_{n+k}}{u_n}\underset{n\to\infty,\ n\notin K_0}\lra\ 1\ \forall\ k\in\mathbb Z$.

To see (\dsagricultural) for $r=s\ \&\ t\in S,\ \ell\in\mathbb Z$,
$$p_{s,t}^{(n+\ell)}=\sum_{k=0}^{n+\ell-1}u_{n+\ell-k}\, {_sp_{s,t}^{(k)}}\ge \sum_{k=0}^{N-1}u_{n+\ell-k}\, {_sp_{s,t}^{(k)}}\ \forall\ n+\ell>N\ge 1.$$
Thus
$$\frac{p_{s,t}^{(n+\ell)}}{u_n}\ge \sum_{k=0}^{N-1}{_sp_{s,t}^{(k)}}\frac{u_{n+\ell-k}}{u_n}\underset{n\to\infty,\ n\notin K_0}\lra\ \sum_{k=0}^{N-1}{_sp_{s,t}^{(k)}}$$
and (\dsagricultural) holds with $K_{s,t}=K_0$. As mentioned above we now have (\Bouquet) with $r=s$.
\

To see (\dsagricultural) for general $r,\ t\in S$, fix first $K_t$ $(u,\mathfrak K)$-small $\st$  $$\frac{p_{s,t}^{(n+k)}}{u_n}\underset{n\to\infty,\ n\notin K_t}\lra\ \pi_t\ \forall\ k\in\mathbb Z.$$
Next,
$$p_{r,t}^{(n+\ell)}=\sum_{k=0}^{n-1} {_sp_{r,s}^{(k)}} p_{s,t}^{(n+\ell-k)}\ge \sum_{k=0}^{N-1} {_sp_{r,s}^{(k)}} p_{s,t}^{(n+\ell-k)}\ \forall\ n+\ell>N\ge 1,$$
and
$$\frac{p_{r,t}^{(n+\ell)}}{u_n}\ge \sum_{k=0}^{N-1}{_sp_{r,s}^{(k)}}\frac{p_{s,t}^{(n+\ell-k)}}{u_n}\underset{n\to\infty,\ n\notin K_t}\lra\ \sum_{k=0}^{N-1}{_sp_{r,s}^{(k)}}\cdot \pi_t\underset{N\to\infty}\lra\  \pi_t$$
($\because\ \ \sum_{k=0}^\infty{_sp_{r,s}^{(k)}}=1$) and (\dsagricultural) holds with $K_{r,t}=K_t$.\ \ \ \CheckedBox\ (\Bouquet).
\

\

\pf{\it of theorem 8.1}:
\

Suppose that $A=[a_1,\dots,a_I]_k\ \&\ B=[b_1,\dots,b_J]_\ell\in\mathcal C_\a$, then for
$n\in\mathbb Z$,
$$A\cap T^{-n}B=\{x\in S^\mathbb Z:\ x_{k+i}=a_i\ \forall\ 1\le i\le I\ \&\ x_{n+\ell+j}=b_j\ \forall\ 1\le j\le J\}$$ and for $n>k+I-\ell$,

\begin{align*}m(A\cap T^{-n}B)&=
\pi_{a_1}p_{a_1,a_2}\dots p_{a_{I-1},a_I}p_{a_I,b_1}^{(n+\ell-I)}p_{b_1,b_2}\dots p_{b_{J-1},b_J}\\ &=m(A)m(B)\cdot\frac{p_{a_I,b_1}^{(n+\ell-I)}}{\pi_{b_1}}\end{align*}
whence
\begin{align*}\frac{m(A\cap T^{-n}B)}{u_n}=\frac{m(A)m(B)}{\pi_{b_1}}\cdot\frac{p_{a_I,b_1}^{(n+\ell-I)}}{u_n}
\overset{(u,\mathfrak K)-\text{\tt\tiny d.}}{\underset{n\to\infty}\lra}\  m(A)m(B).\end{align*}
Rational weak mixing follows from lemma C.\ \  \CheckedBox

\subsection*{Smoothness of  renewal sequences}
\

\subsection*{Remark 8.1}
\

 If $u$ is a recurrent, aperiodic renewal sequence, whose  associated  lifetime distribution $f\in\mathcal P(\mathbb N)$  has tails $f([n,\infty))$ which are $(-\g)$-regularly varying with $\g\in (0,1]$, then

$$\frac{nu_n}{a_n}\ \ \ \ \ \begin{cases}&\ \ \underset{n\to\infty}\lra 1\ \ \ \ \ \g=1\ \ \text{\rm by [E]};\\ &\\ &\ \ \ \ \underset{n\to\infty}\lra \g\ \ \ \ \ \tfrac12<\g<1\ \ \text{\rm by [GL]};\\ &\\ &\ \ \overset{\text{\tt\tiny density}}{\underset{n\to\infty}\lra}\ \ \g \ \ \ \ \g\le\tfrac12\ \ \text{\rm by [GL]}.\end{cases}$$
By proposition 3.1, the convergence in the third case (which follows from Lemma 9.2 below) is also in $u$-density.  In all cases,  $u$
 is smooth and any corresponding Markov chain  is rationally weakly mixing by theorem 8.1.

\proclaim{Proposition 8.3}
\

Suppose that  $u=(u_0,u_1,\dots)$ is an aperiodic, recurrent, renewal  sequence with lifetime distribution $f\in\mathcal P(\mathbb N)$. Let  $L(n):=\sum_{k=1}^nf([k,\infty))$ and $V(t):=\sum_{1\le n\le t} n^2f_n$.

\

\sms{\rm(i)}\ If for some $N\ge 1,\ \sum_{n=N}^\infty\frac1{V(n)^2}<\infty,$ then
$\sum_{n=1}^\infty(u_n-u_{n+1})^2<\infty$.
\

\sms{\rm(ii)}\ If, in addition,  $\frac{L(n)}{\sqrt n}\ \ \underset{n\to\infty}\lra\ \ 0$,
then $u$ is smooth.
\endproclaim
\pf\ {\it of} (i): By  Parseval's formula, and the renewal equation,
\begin{align*} \int_{-\pi}^\pi\frac{|\th|^2d\th}{|1-f(\th)|^2}<\infty\ \ \iff\ \sum_{n=1}^\infty(u_n-u_{n+1})^2<\infty\end{align*}
where $f(\th):=\sum_{n=1}^\infty f_ne^{in\th}$. By aperiodicity, $\sup_{\e\le|\th|\le\pi}|f(\th)|<1\ \forall\ \e>0$ whence
 (using symmetry)
 \begin{align*}\int_{-\pi}^\pi\frac{|\th|^2d\th}{|1-f(\th)|^2}<\infty\ \ \iff\ \ \int_0^{\e}\frac{\th^2d\th}{|1-f(\th)|^2}<\infty\ \text{\tt\small for some}\ \ \e>0.\end{align*}
 Next,
 \begin{align*}|1-f(\th)|&\ge\text{\tt Re}\,(1-f(\th))=
 2\sum_{n=1}^\infty f_n\sin^2(\frac{n\th}2)\\ &\ge
 2\sum_{1\le n\le \frac{\pi}\th} f_n\sin^2(\frac{n\th}2)\ge\frac{2\th^2}{\pi^2}\sum_{1\le n\le \frac{\pi}\th} n^2f_n=:C\th^2V(\frac{\pi}\th).\end{align*}
\

For large $N,\ V(N)>0$ and
\begin{align*}\int_0^{\frac{\pi}N}\frac{\th^2d\th}{|1-f(\th)|^2}\le
\sum_{n=N}^\infty\int_{\frac{\pi}{n+1}}^{\frac{\pi}n}\frac{d\th}{(C\th V(\frac{\pi}\th))^2}\le C' \sum_{n=N}^\infty \frac1{V(n)^2}.\ \ \ \CheckedBox\text{\rm (i)}\end{align*}
\

\pf\ {\it of} (ii):
\

It follows from the renewal equation (see lemma 3.8.5 of [A]) that $a_u(n)\asymp\frac{n}{L(n)}$ whence
$$\frac{\sqrt n}{a_u(n)}\asymp\sqrt n\cdot\frac{L(n)}{n}=\frac{L(n)}{\sqrt n}\underset{n\to\infty}\lra\ 0$$
whence, by (i)
$$\frac1{a_u(n)}\sum_{k=1}^n|u_k-u_{k+1}|\ \le\  \frac{\sqrt n}{a_u(n)}\sqrt{\sum_{n\ge 1}|u_n-u_{n+1}|^2}\ \ {\underset{n\to\infty}\lra}\ \ \ 0.\ \ \CheckedBox\text{\rm (ii)}$$
\

  For example, let  $f\in\mathcal P(\mathbb N)$ be the winnings distribution in the St Petersburg game:
$$f_{k}=\begin{cases}&\tfrac1{2^{n+1}}\ \ \ \ \ \ k=2^{n}\ \ (n\ge 0),\\ &0\ \ \ \ \ \ \text{\tt else.} \end{cases}$$
 The  associated aperiodic, recurrent renewal sequence is  smooth by proposition 8.3 (remark 8.1 above does not apply).

\

The following is ``extends" Dyson's example (on  p. 55 of [Ch]) of an aperiodic renewal sequence  without the strong ratio limit property:
\proclaim{Proposition 8.4}
\

There is  a subsequence smooth, recurrent, renewal sequence which  does not have the strong ratio limit property.
 \endproclaim\demo{Proof} We consider
$\mathbb P:=\{f\in\mathcal P(\mathbb N):\ f_1>0\}$ metrized by $$d(f,g):=|\tfrac1{f_1}-\tfrac1{g_1}|+\sum_{n\ge 1}|f_n-g_n|.$$
This space  is Polish (complete and separable).
\

For $f\in\mathbb P$, let $u^{(f)}$ be the associated (aperiodic, recurrent) renewal sequence. Let
$$\mathbb P_{\text{\tt\tiny SRLP}}:=\{f\in\mathbb P:\ u^{(f)}\ \ \text{\tt has the strong ratio limit property}\}$$
and
$$\mathbb P_{\text{\tt\tiny ss}}:=\{f\in\mathbb P:\ \exists\ \mathfrak K\subset\mathbb N,\ u^{(f)}\ \ \text{\tt is $\mathfrak K$-smooth}\}.$$
We show that
\begin{align*}\tag{\dsmathematical}\mathbb P_{\text{\tt\tiny ss}}\setminus\mathbb P_{\text{\tt\tiny SRLP}}\ \ \text{\tt is residual in}\ \ \ \mathbb P\end{align*}
(and therefore not empty).
\

By Baire's theorem, it suffices to show residuality of  $\mathbb P_{\text{\tt\tiny ss}}$ and $\mathbb P\setminus\mathbb P_{\text{\tt\tiny SRLP}}$.

\demo{Proof that $\mathbb P_{\text{\tt\tiny ss}}$ is residual}
For each $n\ge 1$, the function $f\mapsto u_n^{(f)}$ is continuous ($\mathbb P\to\mathbb R$), being
 a polynomial function of $(f_1,f_2,\dots,f_n)$. Thus
$$\mathbb P_{\text{\tt\tiny ss}}=\bigcap_{k=1}^\infty\bigcup_{N=k}^\infty\{f\in\mathbb P:\ \sum_{j=1}^N|u^{(f)}_j-u^{(f)}_{j+1}|<\frac1k\sum_{j=1}^Nu^{(f)}_j\}$$
 is  a $G_\d$ set. By the renewal theorem
$\mathbb P_{\text{\tt\tiny ss}}\supset\mathbb P_{+}:=\{f\in\mathbb P:\ \sum_{n\ge 1}nf_n<\infty\}$ which is dense in $\mathbb P$.\ \ \CheckedBox

\demo{Proof that $\mathbb P\setminus\mathbb P_{\text{\tt\tiny SRLP}}$ is residual} \ \ For each $k\in\mathbb N$,
$$\Pi_k:=\{f\in\mathbb P:\ \exists\ N>k\ \st\ u_{N-1}^{(f)}<\frac{u_{N}^{(f)}}k\}$$ is open.
\

Since  $\mathbb P\setminus\mathbb P_{\text{\tt\tiny SRLP}}\supseteq\bigcap_{k\ge 1}\Pi_k$, it suffices to prove that each $\Pi_k$ is dense.
To this end, fix $k\ge 1,\ f\in\mathbb P\ \&\ \e>0$. We'll show that $\exists\ g\in\Pi_k,\ d(f,g)<2\e$. To this end note first that $\exists\ h\in\mathbb P$ so that
$d(f,h)<\e$ and so that the set $\{n\in\mathbb N:\ h_n>0\}$ is infinite. Using this, find $\ell>k$ so that
$$0<1-H:=\sum_{j=\ell+1}^\infty h_j<\e.$$
For $L> \ell$ define $g^{(L)}\in\mathbb P$ by
$$g^{(L)}_n:=\begin{cases}& h_n\ \ \ \ \ \ \ n\le\ell;\\ & \sum_{j=\ell+1}^\infty h_j\ \ \ \ \ \ \ n=L;
\\ & 0\ \ \ \ \ \ \ \text{\tt else}.\end{cases}$$
We claim that $\forall\ L$ large,
$u^{(g^{(L)})}_{L-1}<\tfrac{1-H}k\le\tfrac{u^{(g^{(L)})}_{L}}k $, whence $g^{(L)}\in\Pi_k$.
\

To see this define the defective renewal sequence $v$ by
$$v_0=1,\ v_n:=\sum_{k=1}^{n\wedge\ell}h_kv_{n-k},$$
then $u^{(g^{(L)})}_j=v_j\ \forall\ 1\le j<L$.
\

Let
$V_r:=\max_{\nu\ge r\ell+1}v_{\nu}.$
For $j\ge 1$,
$$v_{r\ell+j}=\sum_{i=1}^\ell h_iv_{r\ell+j-i}=\sum_{i=1}^\ell h_iv_{(r-1)\ell+j+\ell-i}\le HV_{r-1}$$
whence $V_r\le HV_{r-1}$ and
$v_{r\ell}\le H^{r-1}v_\ell$.
Now fix $L>\ell$ so that $v_{L-1}<\tfrac{1-H}k.$
It follows as above that $g^{(L)}\in\Pi_k\ \&\ d(f,g^{(L)})<2\e.$\ \ \CheckedBox

\section*{\S9 Examples with Local limit sets}
\

In this section, we prove a generalization of part of theorem 1.1 in [GL] thereby establishing sufficient conditions for rational weak mixing. It is necessary to deal with essentially non-invertible transformations. By corollary E, rational weak mixing passes to the natural extensions of these non-invertible transformations.
\

Suppose that $(X,\B,m,T)$ is a pointwise dual ergodic, measure preserving transformation (as in [A]) with $\g$-regularly varying return sequence $a(n)=a_n(T)$ ($0<\g<1$). As shown in [A] (chapter 3), $T$ is rationally ergodic, and $T$ is not invertible.
\

By the Darling Kac theorem ([DK], see also chapter 3 in [A])
$$\tfrac1{a(n)}S_n(f)\overset{\mathfrak d}\lra \ X_\g m(f)\ \forall\ f\in L^1_+$$ on $(X,\B,m)$ where
$X_\g$ is the Mittag-Leffler distribution of order $\g$ normalized so that $E(X_\g)=1$,  $m(f):=\int_Xfdm$ and $\overset{\mathfrak d}\lra$ on $(X,\B,m)$
denotes convergence in distribution with respect to all $m$-absolutely continuous probabilities.

Let $\Om\in\B,\ m(\Om)=1$ (the normalization $m(\Om)=1$ is not necessary, but convenient).
\

The {\it return time function} to $\Om$ is $\v=\v_\Om:\Om\to\mathbb N$ defined by
$\v(\om):=\min\{n\ge 1:\ T^n\om\in\Om\}<\infty$ a.s. by conservativity. The {\it induced transformation} on  $\Om$ is $T_\Om:\Om\to\Om$ defined by $T_\Om(\om):=T^{\v(\om)}(\om)$. As is well known, $T_\Om$ is an ergodic, probability preserving transformation of $(\Om,\B(\Om),m_\Om)$.
\

The return time process on $\Om$ satisfies the {\tt stable limit theorem}. Indeed, by proposition 1 in [A2],
$$\tfrac1{B(n)}\v_n\overset{\mathfrak d}\lra \ Z_\g $$  on $(\Om,\B(\Om),m_\Om)$ where
$B(n):=a^{-1}(n),\ \ Z_\g=Y_\g^{-\frac1\g}$ is the stable random variable of order $\g$ and $\v_n:=\sum_{k=0}^{n-1}\v\circ T_\Om^k$.
\

The above is true for any $\Om\in\B,\ m(\Om)=1$. By ``choosing" $\Om$ carefully, it may be possible to obtain stronger properties.
\

Accordingly, in the above situation, we call
 $\Om\in R(T),\ m(\Om)=1$  a {\it local limit set} ({\tt LLT}) if $\exists$ a countable, partition  $\b\subset\B(\Om)$ generating $\B(\Om)$ under $T_\Om$ $\st$ $\v_\Om^{-1}\{n\}\in\s(\b)\ \forall\ n\ge 1$ and $\st\ \forall\ \ A,B\in\mathcal C_\b(T_\Om)$,

\begin{align*}\tag{\dsrailways} B(n)m(A\cap T_\Om^{-n}B\cap [\v_n=k_n])\xrightarrow[n\to\infty,\ \frac{k_n}{B(n)}\to x]{}\ f_{Z_\g}(x)m(A)m(B). \end{align*}
uniformly in $x\in [c,d]$ whenever $0<c<d<\infty)$, where $f=f_{Z_\g}$ is the probability density function of $Z_\g$.
\

To be a {\tt LLT} set, essentially, the return time stochastic process to $\Om$ needs to satisfy the conditional,
$\g$-stable, local limit theorem.
\subsection*{Examples 9.1}\ \
\

If $(X,\B,m,T)$ is the tower over the  a Gibbs Markov fibred system (as in [AD]), or an AFU fibred system (as in [ADSZ])
$(\Om,\mathcal A,P,S,\a)$ with $\a$-measurable height function $\v$ satisfying $E(\v\wedge t)$ regularly varying at infinity with index in $(0,1)$, then
$(X,{\mathcal  B},m,T)$ is pointwise dual ergodic, $\Om\in R(T)$ with $a_n(T)=a_n(\Om)\propto\tfrac{n}{E(\v\wedge n)}$ and the return time stochastic process to $\Om$  satisfies the conditional,
$\g$-stable, local limit theorem.  See [AD] and [ADSZ] respectively. Thus, $\Om$  is a {\tt LLT} set.

\

\

\proclaim{Theorem 9.1}
\

Suppose that  $(X,{\mathcal  B},m,T)$ is pointwise dual ergodic with $a(n)=a_n(T)$  $\g$-regularly varying {\rm($\g\in (0,1)$)} and which has a {\tt LLT} set, then  $(X,{\mathcal  B},m,T)$ is  rationally weakly mixing.
\endproclaim
\pf
\

Let $\Om\in R(T)$ be a {\tt LLT} set with accompanying $T_\Om$-generating partition $\b$.  By standardness, up to isomorphism, $\Om=\a^\mathbb N$,  $T_\Om:\Om\to \Om$
is the shift and the collection $\mathcal C_\b(T_\Om)$ of $(\b,T_\Om)$-cylinder sets forms a base of clopen sets for the Polish topology on $\Om$.
 The proof is via lemma C, whose use is enabled by
 the following lemma 9.2, which  is a version of the ``{\tt local limit}" proof of theorem 1.1 of [GL]. Analogous results are established   in [MT].
\

\proclaim{Lemma 9.2}
\

Suppose that  $(X,{\mathcal  B},m,T)$ is pointwise dual ergodic with return sequence   $a(n)=a_n(T)$  which is $\g$-regularly varying {\rm($\g\in (0,1)$)} and which has a {\tt LLT} set
 $\Om\in R(T),\ m(\Om)=1$, then
\begin{align*}\tag{GL}\varliminf_{n\to\infty}\frac{m(A\cap T^{-n}B)}{u_n}\ge m(A)m(B)\ \ \forall\ A,\ B\in\mathcal C_\b(T_\Om)\end{align*} where $u_n:=\tfrac{\g a(n)}n$ and $\b$ is the accompanying $T_\Om$-generating partition.\endproclaim
 \pf (as in [GL]):
 \

 \ \  Fix $A,B\in\mathcal C_\b(T_\Om)$ and $0<c<d<\infty$. Writing $x_{k,n}:=\tfrac{n}{B(k)}$ for $1\le k\le n$ and using the {\tt LLT} property of $\Om$, we have,
\begin{align*}m(A\cap T^{-n}B)&=\sum_{k=1}^nm(A\cap T_\Om^{-k}B\cap[\v_k=n])\\ &\ge
\sum_{1\le k\le n,\ x_{k,n}\in [c,d]}m(A\cap T_\Om^{-k}B\cap[\v_k=x_{k,n}{B(k)}])\\ &\sim
\sum_{1\le k\le n,\ x_{k,n}\in [c,d]}\tfrac{f(x_{k,n})}{B(k)}m(A)m(B)
\end{align*}
as $n\to\infty$ since $\Om$ is a {\tt LLT} set. We are going to show that the last sum is in fact a Riemann sum.
\

Now,
$$x_{k,n}-x_{k+1,n}=\frac{n}{B(k)}-\frac{n}{B(k+1)}\sim \frac{n}{\g kB(k)}$$
as $k,\ n\to\infty$, $x_{k,n}\in [c,d]$ since $B=a^{-1}$ is $\frac1\g$-regularly varying.
\

Also
$$a(n)=a(x_{k,n}B(k))\sim x_{k,n}^\g a(B(k))\sim x_{k,n}^\g k$$
as $k,\ n\to\infty,\ x_{k,n}\in [c,d]$ by the uniform convergence theorem for regularly varying functions.
so
$$\frac1{B(k)}\sim\frac{\g k}n\cdot (x_{k,n}-x_{k+1,n})\sim \frac{\g a(n)}n\cdot\frac{x_{k,n}-x_{k+1,n}}{x_{k,n}^\g}$$
whence, as $n\to\infty$,
\begin{align*}\sum_{1\le k\le n,\ x_{k,n}\in [c,d]}&\frac{f(x_{k,n})}{B(k)}\sim
\frac{\g a(n)}{n}\sum_{1\le k\le n,\ x_{k,n}\in (c,d)}\frac{(x_{k,n}-x_{k+1,n})}{x_{k,n}^\g}f(x_{k,n})\\ &\sim
\frac{\g a(n)}{n}\int_{[c,d]}\frac{f(x)dx}{x^\g}=\frac{\g a(n)}{n}\mathbb E(1_{[c,d]}(Z_\g)Z_\g^{-\g}).\end{align*}
Now $$\mathbb E(1_{[c,d]}(Z_\g)Z_\g^{-\g})=\mathbb E(1_{[c,d]}(X_\g^{-1/\g})X_\g )\underset{c\to 0+,\ d\to\infty}\lra \mathbb E(X_\g )=1,$$
$$\therefore\ m(A\cap T^{-n}B)\gtrsim\ \ \frac{\g a(n)}{n}m(A)m(B).\ \ \CheckedBox$$
Theorem 9.1 now follows from lemma C.\ \ \CheckedBox
\

\

\subsection*{Remark 9.1} \ \ \
\

In some cases,    $\varliminf$ in lemma 9.2 is actually $\lim$ and the transformation has Krickeberg's mixing property. This occurs in:
\sms (i)  the Markov case when $\g\in (\tfrac12,1]$ (in remark 8.1), see  [GL] for $\g\in (\tfrac12,1)$ and [E] for $\g=1$ (see also [Fre]);
\sms (ii) examples  9.1 when $\g\in (\tfrac12,1)$ and sometimes when $\g=1$ (in theorem 9.1), see [MT].
\section*{\S10 Proof of theorem F}
\demo{Proof of (i)}\ \ Recall from [A2] that for $T\in\text{\tt MPT},\ d_k>0\ \&\ Y$ a random variable on $[0,\infty]$,
$\frac{S^{(T)}_{n_k}}{d_k}\overset{\mathfrak d}\lra\ Y$ if
$$g(\tfrac{S^{(T)}_{n_k}(f)}{d_k})\underset{k\to\infty}\lra E(g(Y\int_Xfdm))\ \ \text{\tt weak-* in}\ L^\infty\ \ \ \forall\ g\in C([0,\infty]).$$
The sequence $\{m_j\}$ is called a {\it loose sequence} for $T$ if
$$n_k=m_{j_k}\to\infty,\ d_k>0,\ \frac{S^{(T)}_{n_k}}{d_k}\overset{\mathfrak d}\lra\ Y\ \Lra\ \text{\tt Prob}([Y\in (0,\infty)])=0.$$
As shown  in the proof of theorem 2 in [A2], the collection
$$\{T\in\text{\tt MPT}:\ T\ \ \text{\tt has a loose sequence}\}$$
is residual in {\tt MPT}. No weakly rationally ergodic transformation has a loose sequence and so the collection of these
is meagre in {\tt MPT}. Thus {\tt RWM} is contained in a meagre collection.\ \ \CheckedBox

\

We commence the proof of (ii) by showing:
\demo{Subsequence, rational, weak mixing is residual}
\

We'll use the
\proclaim{Conjugacy Lemma \ \ {\rm (see e.g. [A], [Kri2], [S])}}
\

For  aperiodic $T\in\text{\tt MPT}$,
$$\{\psi^{-1}\circ T\circ \psi:\ \psi\in\text{\tt MPT}\}$$ is dense in {\tt MPT}. \endproclaim

\

By the isomorphism theorem, we may assume WLOG that $(X,\B,m)$ is as in Hopf's example:
$$X=\mathbb R_+\x [0,1],\ \B=\B(\mathbb R_+\x [0,1])\ \&\ \ m=\text{\tt Leb.}.$$
A {\it dyadic square} in $X$ is a square $S=I\x J$ with $I,\ J$ dyadic intervals in $\mathbb R$ (i.e. $\bdy I,\ \bdy J\in\mathbb Q_2$) of the same length.
 A {\it dyadic set} in $X$ is a finite union of dyadic squares. Let $\mathcal D:=\{\text{\tt dyadic sets in}\ X\}$.
\

We'll need the (standard) result that for $N\ge 2$ there is a measure space isomorphism $\Phi_N:X^N\to X$ so that
$$\Phi_N^{-1}(\mathcal D)=\{\text{\tt finite unions of sets in}\ \underset{{N-\text{\tt\tiny times}}}{\underbrace{\mathcal D\otimes\dots\otimes\mathcal D}}\}.$$
\subsection*{Permutations}
\

An automorphism  $T\in\text{\tt MPT}(X,\B,m)$ is called a {\it permutation} if there exist finitely
many disjoint dyadic squares $E_1,\dots,E_N$ and a permutation $\s:\{1,\dots,N\}\to\{1,\dots,N\}$ so that
\bul $T$ maps each $E_i$ onto $E_{\s(i)}$;
\bul $T(x)=x\ \forall\ x\notin\bigcup_{i=1}^NE_i$.
\

The proof of  Satz 2 in [Kri2] applies to show that the collection
\

$\Pi:=\{\text{\tt permutations in MPT}\}$  is dense in {\tt MPT}.
This immediately implies the
\proclaim{Permutation Conjugacy Lemma }
\

For  aperiodic $T\in\text{\tt MPT}$,
$$\{\psi^{-1}\circ T\circ \psi:\ \psi\in\Pi\}$$ is dense in {\tt MPT}. \endproclaim
Note that  $\psi\mathcal R_b=\mathcal R_b$ (the collection  of bounded measurable sets) for $\psi\in\Pi$, but not for arbitrary $\psi\in\text{\tt MPT}$.

\subsection*{Markov shifts in {\tt MPT}$(X)$}
\

We show that any conservative, ergodic, stationary Markov shift with infinite stationary distribution is isomorphic to a piecewise affine transformation
$T\in\text{\tt MPT}(X)$ with a Markov partition whose cylinder sets are  bounded rectangles in $X$.
\

We consider (WLOG) only Markov chains with state space $\mathbb N$.
\

Let $P:\Bbb N\x\Bbb N\to [0,1]$ be a stochastic matrix with infinite stationary distribution
$\pi:\Bbb N\to\Bbb R_+$.
\

We show first that the one-sided shift of $(P,\pi)$ is isomorphic to a measure preserving, piecewise affine map $\tau_{(P,\pi)}:\Bbb R_+\to\Bbb R_+$.
To this end, let
\bul $\a=\{a_k:\ k\in\Bbb N\}$ be a partition $\mod 0$ of $\Bbb R_+$ into open intervals so that $\l(a_s)=\pi_s\ \forall\ s\in\Bbb N$ where
$\l$ denotes Lebesgue measure on $\Bbb R_+$; and
\bul for each $s\in\mathbb N$ let $\{a_{s,t}:\ t\in\mathbb N,\ p_{s,t}>0\}$ be a partition $\mod 0$ of $a_s$ into open intervals so that
$\l(a_{s,t})=\pi_sp_{s,t}\ \forall\ t\in\Bbb N$.
\

Now define $\tau:\Bbb R_+\to\Bbb R_+$ by
$$\tau(x):=\frac{\pi_t}{\pi_sp_{s,t}}\cdot x+\g_{s,t}\ \ \ \ \ \ \ \ x\in a_{s,t}\ \ \ \ \ \ (s,t\in\mathbb N,\ p_{s,t}>0)$$
where $\g_{s,t}$ is chosen so that $\tau a_{s,t}=a_t$.
\

It is standard to show that $\tau_{(P,\pi)}$ preserves $\l$ and is isomorphic to the one-sided shift of $(P,\pi)$.
\

To obtain the two-sided shift of $(P,\pi)$, we represent the {\tt natural extension} of $\tau$ on $\Bbb R_+\x [0,1]$.
\

For $s\in\Bbb N,$ define $v_{a_s}:\tau a_s\to a_s$ by
$$v_{a_s}(y):=\frac{\pi_sp_{s,t}}{\pi_t}\cdot (y-\g_{s,t})\ \ \ \ \ \ \ y\in a_t\subseteq\tau(a_s).$$
Note that $v_{a_s}'=\sum_{t\in\Bbb N,\ p_{s,t}>0}\frac{\pi_sp_{s,t}}{\pi_t}1_{a_t}$.
\

Define for $x\in \Bbb R_+$
$$q_0(x):=0,\ q_k(x):=\sum_{1\le j\le k}1_{\tau a_j}(x)v_{a_j}'(x)\ \ (k\ge 1)$$ and let
$F_{x,a_k}:[0,1]\to [q_{k-1}(x),q_k(x)]$ be the increasing affine map
$$F_{x,a_k}(y):=1_{\tau a_k}(x)v_{a_k}'(x)y+q_{k-1}(x).$$
Now define $T=T_{(P,\pi)}:\Bbb R_+\x [0,1]\to\Bbb R_+\x [0,1]$ by
$$T(x,y):=(\tau (x),F_{\tau x,\a(x)}(y))\ \ \text{\tt\small where}\ \ x\in\a(x)\in\a.$$
 It is standard to show that  $T_{(P,\pi)}\in\text{\tt MPT}\,(X)$   is a natural extension of $\tau$, whence isomorphic to the two-sided shift of $(P,\pi)$.
The partition $\b:=\a\x [0,1]$ is a Markov partition whose cylinder sets are finite unions of bounded rectangles whence $\mathcal{HR(C}_\b)=\mathcal R_b$.

\

\par Let  $P:S\x S\to [0,1]$ be a  stochastic matrix on the  state space $S$  with  invariant distribution $\pi:S\to \mathbb R_+$.
and let $T_{(P,\pi)}$  be a Markov shift in {\tt MPT} isomorphic to the
 stationary  Markov shift of $(P,\pi)$. Fix $s\in S$ and let $u=u([s]_0),\ a(n):=\sum_{k=0}^{n-1}u_k$.

\

Assume that $u$ is smooth, then $T=T_{(P,\pi)}$ is rationally weakly mixing with $R(T)\supset\mathcal{HR(C}_\b)=\mathcal R_b$, whence

\begin{align*}\tag{$\maltese$}\label{maltese1}&\frac1{a(n)}\sum_{k=0}^{n-1}|m(D\cap T^{-n}D')-u_nm(D)m(D')|\underset{n\to\infty}\lra\ 0\ \forall\ D,\ D'\in\mathcal D;\end{align*}
which implies
\begin{align*}\frac1{a(n)}\int_DS_n^{(T)}(1_D)dm\underset{n\to\infty}\lra m(D)^2\ \ \ \forall\ D\in\mathcal D.\end{align*}
We claim that also
\begin{align*}\tag{{\Large\Bat}}\label{bat} \varlimsup_{n\to\infty}\frac1{a(n)^2}\int_DS_n^{(T)}(1_D)^2dm\ \le 2m(D)^3\ \ \ \forall\ D\in\mathcal D.\end{align*}
\demo{Proof of ({\Large\Bat})}\ \ \ Let $\b=\a\x [0,1]$, the Markov partition of $T$.
We first show  ({\Large\Bat}) for $A\in\mathcal U_\b$. Let $\tau=\tau_{(P,\pi)}:\Bbb R_+\to\Bbb R_+$ be as above (isomorphic to the
  one-sided Markov shift  of $(P,\pi)$). It is pointwise dual ergodic in the sense that
 $$\frac1{a(n)}\sum_{k=0}^{n-1}\widehat{\tau} ^k1_A\underset{n\to\infty}\lra\ \l(A)\ \forall\ A\in\mathcal F\tag{a}$$ where
 $\widehat{\tau} :L^1(\l)\to L^1(\l)$ is the transfer operator defined by $$\int_{\Bbb R_+}\widehat{\tau} f\cdot gd\l= \int_{\Bbb R_+}f\cdot g\circ\tau d\l$$ (see \S3.7 in [A]).
 \

Now $\sup_{\Bbb R_+}\sum_{k=0}^{n-1}\widehat{\tau} ^k1_A=\sup_{A}\sum_{k=0}^{n-1}\widehat{\tau} ^k1_A$.
 For $A\in\mathcal C_{\a}$, the convergence (a) is uniform on $A$, whence
 $$\frac1{a(n)}\sup_{\Bbb R_+}\sum_{k=0}^{n-1}\widehat{\tau} ^k1_A=\frac1{a(n)}\sup_{A}\sum_{k=0}^{n-1}\widehat{\tau} ^k1_A\underset{n\to\infty}\lra\ m(A)\ \forall\ A\in\mathcal C_{\a}.\tag{b}$$
 From (b) we see that
$$\varlimsup_{n\to\infty}\frac1{a(n)}\sup_{\Bbb R_+}\sum_{k=0}^{n-1}\widehat{\tau} ^k1_A\ \le\ m(A)\ \forall\ A\in\mathcal C_{\a}.\tag{b'}$$

The statement (b$^\prime$)  holds  $$\forall\ A\in \mathcal U_{\a}:=\{\bigcup_{k=1}^NC_k:\ \ N\ge 1,\ C_1,\dots,C_N\in\mathcal C_{\a}\}$$ and it follows that for $A\in \mathcal U_{\a}$, and $n$ large so that
 $\sum_{k=0}^{n-1}\widehat{\tau} ^k1_A\le 2m(A)a(n)$,
 \begin{align*}\tag{c}\int_AS_n(1_A)^2d\l &\le 2\sum_{0\le i\le j\le n-1}\l (A\cap \tau^ {-i}A\cap \tau^ {-j}A)\\ &=
 2\int_A\sum_{0\le i\le n-1}\widehat{\tau} ^i1_AS_{n-i}(1_A)\circ \tau^ id\l \\ &\le\int_AS_n(1_A)\sum_{k=0}^{n-1}\widehat{\tau} ^k1_Ad\l \\ &\le
 4m(A)a(n)\int_AS_n(1_A)d\l \\ &\le 8m(A)^3a(n)^2.\end{align*}
If $\psi:X=\Bbb R_+\x [0,1]\to \Bbb R_+$ is the projection $\psi(x,y)=x$, then
$$\mathcal U_\b=\bigcup_{n\ge 1}T^n\psi^{-1}\mathcal U_{\a}$$

and  ({\Large\Bat}) follows for $A\in\mathcal U_\b$. Dyadic sets can be monotonically approximated by sets in $\mathcal U_\b$ and  ({\Large\Bat}) follows.\ \ \CheckedBox

\

Now enumerate $\mathcal D:=\{D_i:\ i\in\mathbb N\}$ and define
\begin{align*}\mathcal G:=
\bigcap_{k\ge 1}\bigcup_{N\ge k}\bigcap_{1\le i,j\le k}U(N,k,i,j)\end{align*}
where
\begin{align*}U(N,k,&i,j):=\\ &\{T\in\text{\tt MPT}:\ \sum_{\nu=0}^{N-1}|m(D_i\cap T^{-\nu}D_j)-m(D_i)m(D_j)u_\nu|<\frac{a(N)}k;\\ &\ \ \ \ \ \ \ \ \ \ \ \ \  \& \ \ \int_{D_i}S_N(1_{D_i})^2dm<8m(D_i)^3a(N)^2\}.\end{align*}
Evidently each $U(N,k,i,j)$ is open in {\tt MPT}, whence the set $\mathcal G$ is a $G_\d$ set. We'll complete the proof of residuality of {\tt SRWM}  by showing that
\bul $\mathcal G$ is dense in {\tt MPT} and
\bul each $T\in\mathcal G$ is subsequence rationally weakly mixing.
\demo{Proof of density of $\mathcal G$}

By ($\maltese$)  and ({\Large\Bat}) (as on page  \pageref{maltese1}), $T_{(P,\pi)}\in\mathcal G$. Since $\psi\mathcal R_b=\mathcal R_b\ \forall\ \psi\in\Pi$,
$$\{\psi^{-1}\circ T_{(P,\pi)}\circ \psi:\ \psi\in\Pi\}\subset\mathcal G.$$
Since $T_{(P,\pi)}$ is ergodic, by the  permutation conjugacy lemma,
\begin{align*}\o{\mathcal G}\supset\o{\{\psi^{-1}\circ T_{(P,\pi)}\circ \psi:\ \psi\in\Pi\}}=\text{\tt MPT}.\ \ \CheckedBox\end{align*}
\demo{Proof of subsequence rational weak mixing of elements of $\mathcal G$}
Let $T\in\mathcal G$, then $\exists$ a subsequence $\mathfrak K\subset\mathbb N$ such that
\begin{align*}&\tag{a}
 \tfrac1{a(N)}\sum_{\nu=0}^{N-1}|m(D\cap T^{-\nu}D')-m(D)m(D')u_\nu|\underset{N\to\infty,\ N\in\mathfrak K}\lra 0\ \forall\ D,D'\in\mathcal D;
\\ &\tag{b}\ \ \  \int_{D}S_N(1_{D})^2dm<8m(D)^3a(N)^2\ \ \ \forall\ D\in\mathcal D,\ N\in\mathfrak K.\end{align*}
It follows from (a) that
$$\frac{a_N^{(T)}({D})}{a(N)}\underset{N\to\infty,\ N\in\mathfrak K}\lra m(D)^2\ \ \ \forall\ D\in\mathcal D$$
whence by
(b), $T$ is weakly rationally ergodic along $\mathfrak K$ with return sequence $a(n)$ along $\mathfrak K$
 and $\mathcal D\subset R_{\mathfrak K}(T)$. This enables use of  (a) and lemma C to show that $T$ is rationally weakly mixing along $\mathfrak K$.\ \ \CheckedBox

\demo{Proof of (ii)}
\

For  $\u\kappa=(\kappa_1,\dots,\kappa_\D)\in(\Bbb Z\setminus\{0\})^\D,$ define $\psi_{\u\kappa}:\text{\tt MPT}\,(X)\to\text{\tt MPT}\,(X)$ by
$$\psi_{\u\kappa}(T):=\phi_\D\circ T^{(\kappa_1,\dots,\kappa_\D)}\circ\phi_\D^{-1}\in\text{\tt MPT}(X)$$
where as above, $\Phi_{\D}:X^{\D}\to X$ so that
$$\Phi_{\D}^{-1}(\mathcal D)=\{\text{\tt finite unions of sets in}\ \underset{\D-\text{\tt\tiny times}}{\underbrace{\mathcal D\otimes\dots\otimes\mathcal D}}\}.$$
If $\psi_{\u\kappa}(T)\in\mathcal G$, then  $T^{(\kappa_1,\dots,\kappa_\D)}\in\text{\tt SRWM}$. Thus it suffices to show that
$$\mathcal G_{\text{\tt\tiny power}}:=\bigcap_{\D=1}^\infty\ \bigcap_{(\kappa_1,\dots,\kappa_\D)\in(\Bbb Z\setminus\{0\})^\D}\psi_{\u\kappa}^{-1}\mathcal G$$ is residual.

\

It is not hard to see that:\bul   each $\psi_{\u\kappa}:\text{\tt MPT}\,(X)\to\text{\tt MPT}\,(X)$ is a continuous homomorphism, whence $\mathcal G_{\text{\tt\tiny power}}$ is a $G_\d$ set in
{\tt MPT}$\,(X)$;  and that \bul $\psi_{\u\kappa}(\Pi)=\Pi$, whence $\psi^{-1}\circ T\circ\psi\in\mathcal G_{\text{\tt\tiny power}}\
\forall\ T\in\mathcal G_{\text{\tt\tiny power}},\ \psi\in\Pi$, because for
 $T\in\mathcal G_{\text{\tt\tiny power}}\ \&\ \pi\in\Pi$, $\psi_{\u\kappa}(\pi)\mathcal D=\mathcal D$ and
$$\psi_{\u\kappa}(\pi^{-1}\circ T\circ\pi)=\psi_{\u\kappa}(\pi)^{-1}\circ\psi_{\u\kappa}(T)\circ\psi_{\u\kappa}(\pi)\in\mathcal G.$$
\

To prove density of $\mathcal G_{\text{\tt\tiny power}}$ (and thus complete the proof of (ii)) it suffices to exhibit $T\in\mathcal G_{\text{\tt\tiny power}}$ for then
 $T$ is ergodic and
\

$\o{\mathcal G_{\text{\tt\tiny power}}}\supset\o{\{\pi^{-1}\circ T\circ\pi:\ \pi\in\Pi\}}=\text{\tt MPT}$ by the permutation conjugacy lemma.
\

\subsection*{Renewal shifts}\ \
 Let $u$ be a   recurrent, renewal sequence with
 lifetime distribution $f\in\mathcal P(\Bbb N)$.
Define (as in [Ch]) a stochastic matrix $P=P_u$ on $\Bbb N$ by
$$P_{1,n}:=f_n\ \ \&\ \ P_{n+1,n}=1\ \ \ \ \forall\ n\in\Bbb N.$$
This has stationary distribution $\pi=\pi_u$ defined by
$\pi_n:=\sum_{k=n}^\infty f_k$ and $P_{1,1}^{(n)}=u_n$. The
 Markov shift of $(P,\pi)$ is called the {\it renewal shift} of $u$. Let $T_u:=T_{(P,\pi)}\in\text{\tt MPT}$.
\

If $u$ is smooth, then
$T_u\in\text{\tt RWM}$.
\

Now suppose that $\D\ge 1,\ \u\kappa=(\kappa_1,\dots,\kappa_\D)\in(\Bbb Z\setminus\{0\})^\D$, then as evidently $T_u^{-1}\cong T_u$,
$$T_u^{\kappa_1}\x\dots\x T_u^{\kappa_\D}\cong T_u^{|\kappa_1|}\x\dots\x T_u^{|\kappa_\D|}$$ and we may assume WLOG that $\u\kappa\in\Bbb N^\D$. Now
 $T_u^{\kappa_1}\x\dots\x T_u^{\kappa_\D}$ is also the Markov shift of an irreducible, aperiodic, stochastic matrix with renewal sequence $u^{(\u\kappa)}$ defined by
$$u^{(\u\kappa)}_n:=\prod_{j=1}^\D u_{\kappa_jn}.$$
\

If $u$ is smooth and  $u^{(\u\kappa)}$ is recurrent, then  $u^{(\u\kappa)}$ is also smooth,
\

  $T_u^{\kappa_1}\x\dots\x T_u^{\kappa_\D}$ is rationally weakly mixing and   $\psi_{\u\kappa}(T_u)\in\text{\tt RWM}.$

\

Now let $u$ be the sequence defined by
$$u_n:=\frac1{\log(n+e)}\ \ \ \ \ \ (n\ge 0),$$
then $u$ is a {\it Kaluza sequence} in the sense that $u_0=1\ \&\ \tfrac{u_{n+1}}{u_n}\uparrow 1$ and hence a smooth, recurrent renewal sequence.
\

As can be easily checked, so is $u^{(\u\kappa)}$ $\forall\ \D\ge 1,\ \u\kappa=(\kappa_1,\dots,\kappa_\D)\in\Bbb N^\D$. It follows that
$T_u\in\mathcal G_{\text{\tt\tiny power}}$.\ \ \ \CheckedBox

\section*{\S11 Closing Remarks}
\

 All  infinite, rationally weakly mixing  examples in this paper are of form $T\x S$ where $T$ is an infinite K-automorphism and
$S$ is a weakly mixing probability preserving transformation.
\

Their Koopman operators all have countable Lebesgue spectrum.
This is shown in [Par] for   K-automorphisms and a simple argument shows that multiplying by a weakly mixing probability preserving
transformation does not change this.
\

The transformation $T\in\text{\tt MPT}$ is called {\it rigid} if $\exists\ \frak L\subset\Bbb N$ so that
$$m(A\D T^{-n}A)\underset{n\to\infty,\ n\in\frak L}\lra\ 0\ \ \ \forall\ A\in\mathcal F.$$
The spectrum of a rigid transformation is Dirichlet, and hence singular.

As shown in  [AS], the collection $\text{\tt RIGID}$ of rigid transformations in {\tt MPT} is residual.
\

By Theorem F, so is $\text{\tt PSRWM}\cap\text{\tt RIGID}$ and so there is a rigid, power, subsequence, rationally weak mixing, measure preserving transformation with singular spectrum.
 
\end{document}